\documentclass{gtart_a}
\pdfoutput=1

\usepackage[all]{xy}
\usepackage{graphicx}

\title[Euler structures, the variety of representations and the Milnor--Turaev torsion]{Euler structures, the variety of representations\\and the Milnor--Turaev torsion}

\author{Dan Burghelea}
\givenname{Dan}
\surname{Burghelea}
\address{Dept. of Mathematics\\
The Ohio State University\\\newline
231 West 18th Avenue\\
Columbus, OH 43210\\USA}
\email{burghele@mps.ohio-state.edu}
\urladdr{}

\author{Stefan Haller}
\givenname{Stefan}
\surname{Haller}
\address{Department of Mathematics\\ 
University of Vienna\\\newline
Nordbergstrasse 15\\
A-1090, Vienna\\Austria}
\email{stefan.haller@univie.ac.at}
\urladdr{}

\volumenumber{10}
\issuenumber{}
\publicationyear{2006}
\papernumber{29}
\lognumber{0689}
\startpage{1185}
\endpage{1238}

\doi{}
\MR{}
\Zbl{}

\arxivreference{math.DG/0310154}


\received{16 December 2005}
\revised{27 March 2006}
\accepted{23 July 2006}
\published{16 September 2006}
\publishedonline{16 September 2006}
\proposed{Wolfgang L\"uck}
\seconded{Bill Dwyer, Haynes Miller}
\corresponding{}
\editor{}
\version{}

\let\xysavmatrix\xymatrix
\def\xymatrix{\disablesubscriptcorrection\xysavmatrix}
\AtBeginDocument{\let\bar\wbar\let\tilde\wtilde\let\hat\what}
\newcommand\RS{\text{\rm {\tiny RS}}}

\makeatletter
\def\cnewtheorem#1[#2]#3{\newtheorem{#1}{#3}[section]
\expandafter\let\csname c@#1\endcsname\c@theorem}

\newcommand\itemref[1]{\eqref{#1}}
\renewcommand{\theenumi}{\roman{enumi}}

\theoremstyle{plain}
  \newtheorem{theorem}{Theorem}[section]
  \cnewtheorem{corollary}[theorem]{Corollary}
  \cnewtheorem{lemma}[theorem]{Lemma}
  \cnewtheorem{proposition}[theorem]{Proposition}

\theoremstyle{definition}
  \cnewtheorem{definition}[theorem]{Definition}
  \cnewtheorem{remark}[theorem]{Remark}
  \cnewtheorem{example}[theorem]{Example}

\makeatother  

\newcommand{\an}{\text{\rm an}}
\newcommand\even{\text{\rm even}}
\newcommand\odd{\text{\rm odd}}
\newcommand\can{\text{\rm can}}
\newcommand\ac{\text{\rm ac}}
\newcommand\Hodge{\text{\rm Hodge}}
\newcommand\cc{\mathfrak c}
\renewcommand\aa{\mathfrak a}
\newcommand\bb{\mathfrak b}
\newcommand\hh{\mathfrak h}
\newcommand\id{\text{\rm id}}
\newcommand{\N}{\mathbb N}
\newcommand{\X}{\mathbb X}
\newcommand{\T}{\mathbb T}
\newcommand{\K}{\mathbb K}
\newcommand{\Eul}{\mathfrak{Eul}}
\newcommand\vf{\mathfrak X}
\newcommand{\ie}{ie, }
\newcommand{\e}{\mathfrak e}
\newcommand\Or{\mathcal O}
\newcommand\comb{\text{\rm comb}}
\DeclareMathOperator{\rank}{rank}
\DeclareMathOperator{\GL}{GL}
\DeclareMathOperator\Hom{Hom}
\DeclareMathOperator\coker{coker}
\DeclareMathOperator\img{img}
\DeclareMathOperator\Emb{Emb}
\DeclareMathOperator\tr{tr}
\DeclareMathOperator{\Det}{Det}
\DeclareMathOperator{\sign}{sign}
\DeclareMathOperator{\IND}{IND}
\DeclareMathOperator{\End}{End}
\DeclareMathOperator{\cs}{cs}
\DeclareMathOperator\PD{PD}
\DeclareMathOperator\Rep{Rep}

\begin{document}

\begin{asciiabstract}
In this paper we extend and Poincare dualize the concept of Euler
structures, introduced by Turaev for manifolds with vanishing
Euler-Poincare characteristic, to arbitrary manifolds.  We use the
Poincare dual concept, co-Euler structures, to remove all geometric
ambiguities from the Ray-Singer torsion by providing a slightly
modified object which is a topological invariant.  We show that when
the co-Euler structure is integral then the modified Ray-Singer
torsion when regarded as a function on the variety of generically
acyclic complex representations of the fundamental group of the
manifold is the absolute value of a rational function which we call in
this paper the Milnor-Turaev torsion.
\end{asciiabstract}

\begin{abstract}
In this paper we extend and Poincar\'e dualize the concept of Euler 
structures, introduced by Turaev for manifolds with vanishing 
Euler--Poincar\'e characteristic, to arbitrary manifolds.
We use the Poincar\'e dual concept, co-Euler structures, to remove all
geometric ambiguities from the Ray--Singer torsion by providing a slightly
modified object which is a topological invariant.
We show that when the co-Euler structure is integral then the modified 
Ray--Singer torsion when regarded
as a function on the variety of generically acyclic complex representations 
of the fundamental group of the manifold is the absolute value of a rational 
function which we call in this paper the Milnor--Turaev torsion. 
\end{abstract}

\maketitle

\section{Introduction}\label{S:intro}

This paper was motivated by a question in geometric analysis.

{\bf Question}\qua 
{\sl Is the Ray--Singer torsion, when regarded as a positive real-valued function
on the space of complex representations of the fundamental group, 
the absolute value of a holomorphic/meromorphic function. If so,
what is the meaning of its argument?}\footnote{A similar question was 
treated by Quillen \cite{Q85}. Instead of the variety of representation Quillen 
considers the affine space of complex structures in a smooth complex vector 
bundle over a Riemann surface, and instead of the Ray--Singer torsion the square 
root of the determinant of $D^*D$ where $D$ is the appropriate Cauchy--Riemann 
operator. More recently, a series of preprints of Braverman and Kappeler 
\cite{BK06,BK05,BK05b,BK05a} have been posted treating a similar question.}

The answer provided in this paper goes as follows; cf
\fullref{SS:apl}. With the help of an integral co-Euler structure
(for definitions and properties see \fullref{SS:coeul}) we decompose 
the Ray--Singer torsion as a product of two positive real-valued functions, 
the combinatorial torsion and the Bismut--Zhang anomaly.
We verify that the first one, when restricted to generically acyclic representations, 
is the absolute value of a rational function, the \emph{Milnor--Turaev torsion\/}, 
and the second is the absolute value of a holomorphic function.

The rational function is a topological invariant of the manifold and 
carries significant topological information; see \fullref{SS:maptor}. 
Its absolute value is a familiar quantity but the argument is an unexplored 
invariant which is analogous to the Atiyah--Patodi--Singer eta invariant
\cite{APS75}, but not the same.

The values of the second function are obtained from an integral on a 
noncompact manifold, possibly divergent, but always regularizable and 
referred to as the \emph{invariant $R$\/} defined in \fullref{SS:R}.
This invariant was first introduced by Bismut--Zhang and is a sort of 
\emph{Chern--Simons\/} quantity. Both, the regularization of the integral 
and the structure of the invariant $R$, deserve and receive 
in this paper special attention.

The main results of the paper are contained in \fullref{T:eul}, \fullref{T:ratio},
\fullref{T:marcsik}, \fullref{C:meroa} and \fullref{C:merob}. The 
paper contains a number of additional facts/re\-sults of independent
interest expanding on the concepts described in the paper.

Let us make a more detailed presentation of the contents of this paper.

\fullref{char} is largely expository. We collect a number of facts about 
Riemannian metrics and vector fields with isolated zeros and about the 
Kamber--Tondeur form; see \fullref{SS:MQECS}--\fullref{SS:KT}. 
In \fullref{SS:reg} we describe the regularization of some possibly 
divergent integrals which is essential in the subsequent considerations 
(the invariant $R$, the comparison of Euler and co-Euler structures, the 
modifications of torsions). In \fullref{SS:R} we define the invariant $R$.

We regard the invariant $R$ as providing a pairing between a Riemannian metric
and a vector field with isolated zeros, with values in one-currents rel.\ boundary.
This is in analogy with the Chern--Simons form which provides such pairing for 
two Riemannian metrics and with the Chern--Simons one-chain rel.\ boundary which 
provides a pairing between two vector fields with isolated zeros.

\fullref{S:eul} discusses Euler and co-Euler structures and what they are
good for. Turaev \cite{Tu90} introduced the concept of Euler structures for
manifolds with vanishing Euler--Poincar\'e characteristics $\chi(M)$.
It was observed by the first author \cite{B99} that Euler structures can be defined
for manifolds with arbitrary Euler characteristics at the expense of 
a base point $x_0\in M$. The set of Euler structures $\Eul_{x_0}(M;\Z)$
is an affine version of $H_1(M;\Z)$ in the sense that $H_1(M;\Z)$ acts 
freely and transitively on $\Eul_{x_0}(M;\Z)$. There is also a real version 
of Euler structures which is an affine version of $H_1(M;\R)$ and will be 
denoted by $\Eul_{x_0}(M;\R)$. There is a canonical map
$\Eul_{x_0}(M;\Z)\to\Eul_{x_0}(M;\R)$ which is affine over the homomorphism
$H_1(M;\Z)\to H_1(M;\R)$; see \fullref{SS:eul}.

Given an Euler structure $\e\in\Eul_{x_0}(M;\Z)$, a cohomology orientation
$\mathfrak o$, \ie an orientation of $\bigoplus_iH^i(M;\R)$, and a representation 
$\rho\in\Rep(\Gamma;V)$, $\Gamma=\pi_1(M,x_0)$, we use the Milnor--Turaev 
construction \cite{M66,Tu90} to provide an element
$\tau_\comb^{\rho,\e,\mathfrak o}$ in 
the determinant line\footnote{Here we use the
following notation. For a finite-dimensional vector space $V$ we write 
$\det V:=\Lambda^{\dim V}V$ for its determinant line, a one-dimensional vector 
space. For a finite-dimensional graded vector
space $W^*$ we define the (graded) determinant line by $\det W^*:=\det
W^\even\otimes(\det W^\odd)^*$. Moreover, for a line $L$, \ie one-dimensional
vector space, and $k\in\N$ we write $L^k:=L\otimes\cdots\otimes L$
and $L^{-k}:=(L^k)^*$. Then there are natural identifications 
$L^{k+l}=L^k\otimes L^l$ for all $k,l\in\Z$.}
\begin{equation}\label{E:Det}
\Det_{x_0}(M;\rho):=\det H^*(M;\rho)\otimes\det V^{-\chi(M)}
\end{equation}
which will be called \emph{Milnor--Turaev torsion\/}; see \fullref{D:MT-tor}.

For the dependence on $\e$ and $\mathfrak o$ we find 
\begin{equation}\label{E:dep}
\tau^{\rho,\e+\sigma,\mathfrak o}_\comb=\tau^{\rho,\e,\mathfrak o}_\comb\cdot
[\det\circ\rho](\sigma)^{-1}
\quad\text{and}\quad
\tau^{\rho,\e,-\mathfrak o}_\comb=(-1)^{\dim V}\tau^{\rho,\e,\mathfrak o}_\comb
\end{equation}
for all $\sigma\in H_1(M;\Z)$.
Here $[\det\circ\rho]\co H_1(M;\Z)\to\C^*$ denotes the homomorphism induced from
the homomorphism $\det\circ\rho\co\Gamma\to\C^*$.

If $H^*(M;\rho)=0$ then $\Det_{x_0}(M;\rho)=\C$, and
$\tau_\comb^{\rho,\e,\mathfrak o}$ is a nonzero complex number.

In \fullref{SS:coeul}
we introduce the concept of co-Euler structures. The set of co-Euler
structures $\Eul^*_{\smash{x_0}}(M;\R)$ is an affine version of $H^{n-1}(M;\Or_M)$
and depends on the choice of a base point $x_0\in M$ if $\chi(M)\neq0$.
By definition a co-Euler structure is an equivalence class of pairs $(g,\alpha)$, where $g$ 
is a Riemannian metric on $M$, and $\alpha\in\Omega^{n-1}(M\setminus\{x_0\};\Or_M)$ 
with $d\alpha=E(g)$ where $E(g)\in\Omega^n(M;\Or_M)$ denotes the Euler form of $g$. Two 
such pairs $(g_1,\alpha_1)$ and $(g_2,\alpha_2)$ are equivalent if
$\alpha_2-\alpha_1=\cs(g_1,g_2)$ where
$\cs(g_1,g_2)\in\Omega^{n-1}(M;\Or_M)/d\Omega^{n-2}(M;\Or_M)$
denotes the Chern--Simons class. The main properties of Euler and co-Euler
structures and the relationship between them are collected in
\fullref{T:eul}. Co-Euler structures are nicely
suited to remove the metric ambiguities of the analytic torsion and finally
provide a topological invariant, \emph{the modified Ray--Singer metric/torsion\/}.

More precisely, let $\mathbb F_\rho$ denote the flat complex vector bundle 
over $M$ associated to a representation $\rho\in\Rep(\Gamma;V)$.
Recall that the Ray--Singer metric is a Hermitian metric
$\|\cdot\|_\RS^{\rho,g,\mu}$ on the determinant line 
$\det H^*(M;\mathbb F_\rho)=\det H^*(M;\rho)$
depending on a Riemannian metric $g$ on $M$ and a Hermitian structure 
$\mu$ on $\mathbb F_\rho$ when $M$ has even dimension, but is independent  
of them when $M$ is of odd dimension; see Bismut and Zhang \cite{BZ92}. Using a co-Euler structure
$\e^*\in\Eul^*_{\smash{x_0}}(M;\R)$ permits to define a 
\emph{modified Ray--Singer metric\/} $\|\cdot\|_\an^{\smash{\rho,\e^*}}$ on 
the determinant line \eqref{E:Det}; see \fullref{D:mRS-met}.
The anomaly formulas for the Ray--Singer torsion imply that this Hermitian
metric indeed only depends on the representation $\rho$ and the co-Euler
structure $\e^*$. For the dependence on $\e^*$ we have
$$
\|\cdot\|_\an^{\rho,\e^*+\beta}
=\|\cdot\|_\an^{\rho,\e^*}\cdot
e^{\langle\log|[\det\circ\rho]|,\,\PD^{-1}(\beta)\rangle}
$$
for all $\beta\in H^{n-1}(M;\Or_M)$. Here 
$\log|[\det\circ\rho]|\in H^1(M;\R)$ denotes the class
corresponding to the homomorphism $\log|[\det\circ\rho]|\co H_1(M;\Z)\to\R$,
which is paired with the Poincar\'e dual $\PD^{-1}(\beta)$ of $\beta$.

There is a natural isomorphism $P\co\Eul_{x_0}(M;\R)\to\Eul^*_{\smash{x_0}}(M;\R)$ which is
affine over the Poincar\'e duality isomorphism 
$\PD\co H_1(M;\R)\to H^{n-1}(M;\Or_M)$; see \fullref{SS:coeul}. 
We continue to denote by $P$ the composition of $P$ with the canonical map
$\Eul_{x_0}(M;\Z)\to\Eul^*_{\smash{x_0}}(M;\R)$ and refer to $P(\Eul_{x_0}(M;\Z))$ as the 
\emph{lattice of integral co-Euler structures\/}.

This permits a  
reformulation of the Bismut--Zhang theorem \cite{BZ92} as a statement which compares 
two topological invariants; see \fullref{T:BZ}. 
More precisely, it says that analysis is able to recover the absolute value 
of the Milnor--Turaev torsion
$$
\|\tau_\comb^{\rho,\e,\mathfrak o}\|_\an^{\rho,\e^*}
=e^{\langle\log|[\det\circ\rho]|,\,\PD^{-1}(\e^*-P(\e))\rangle}.
$$
Here $\log|[\det\circ\rho]|\in H^1(M;\R)$ denotes the cohomology class
corresponding to the homomorphism $\log|[\det\circ\rho]|\co H_1(M;\Z)\to\R$.

\fullref{S:rep} is dedicated to the Milnor--Turaev torsion viewed as a 
rational function on an algebraic set of generically acyclic representations
which will be described below.
This rational function is presented from a perspective analogous
to the theory of characteristic classes of vector bundles 
(via the Grassmann variety, the Chern forms of the tautological bundle and 
the classifying map.)

Here is a brief description. For a collection $(k_0,k_1,\dotsc,k_n)$
of nonnegative integers satisfying
$k_i-k_{i-1}\pm\cdots+(-1)^ik_0\geq 0$, $i\leq n-1$ and
$k_n-k_{n-1}\pm\cdots+(-1)^nk_0=0$,
we denote by $\mathbb D_\ac(k_0,\dotsc,k_n)$ the set of acyclic
complexes $0\to\mathbb C^{k_0}\to\cdots\to\mathbb C^{k_n}\to 0$.
In \fullref{SS:Tfunc} we will see that its Zariski closure
is an (irreducible) affine variety.
We will denote this closure by $\hat{\mathbb D}_\ac(k_0,\dotsc,k_n)$ and
refer to it as the \emph{variety of generically acyclic complexes\/}.
The assignment which associates to an acyclic complex $C$ its torsion
$\mathfrak t(C)$ defines a rational function $\mathfrak t$ on
$\hat{\mathbb D}_\ac(k_0,\dotsc,k_n)$.

For a finitely presented group $\Gamma$ we denote by $\Rep(\Gamma;V)$
the complex algebraic set of complex representations of $\Gamma$ on
the complex vector space $V$; cf \fullref{SS:rep}.

For a closed pointed manifold $(M,x_0)$ with $\Gamma=\pi_1(M,x_0)$
denote by $\Rep^M(\Gamma;V)$ the algebraic closure of $\Rep^M_0(\Gamma;V)$, that is,
the Zariski open set of representations $\rho\in\Rep(\Gamma;V)$ so that
$H^*(M;\rho)=0$. The elements of $\Rep^M(\Gamma;V)$ are called
\emph{generically acyclic\/} representations.

By choosing a triangulation $\tau$ of $M$, a collection of paths
$\pi_\tau=\{\pi_{x_\sigma}\}$ from $x_0$ to the barycenters $x_\sigma$
of each simplex $\sigma$, an ordering $o$ of the set of simplices
and a base (frame) $\epsilon$ of $V$, one obtains a regular map
$t_{\pi_\tau,o,\epsilon}\co\Rep^M(\Gamma;V)\to
\hat{\mathbb D}_\ac(k_0,\dotsc,k_n)$ with
$k_i=\dim(V)\times\sharp(\mathcal X_i)$ where
$\mathcal X_i$ denotes the set of $i$--simplices and $\sharp$ denotes
cardinality. One can show that while the map $t_{\pi_\tau,o,\epsilon}$
depends on $(\pi_\tau,o,\epsilon)$, the pull-back of $\mathfrak t$ by
$t_{\pi_\tau,o,\epsilon}$
depends only on the Euler structure $\e$ defined by $(\tau,\pi_{\tau})$
(more precisely by $(X_\tau,\pi_\tau)$, cf \fullref{SS:eul}), and the
cohomology orientation $\mathfrak o$ induced by $o$.
We denote the pullback $\mathfrak t\cdot t_{\pi_\tau,o,\epsilon}$ by
$$
\mathcal T_\comb^{\e,\mathfrak o}\co\Rep^M_0(\Gamma;V)\to\C^*.
$$
In \fullref{SS:Tfunc}, in particular \fullref{T:ratio}, we observe that this function is 
a rational function on $\Rep^M(\Gamma;V)$ whose zeros and poles are
therefore contained in the complement $\Rep^M(\Gamma;V)\setminus\Rep^M_o(\Gamma;V)$.
Recall that for acyclic $\rho$ we have a canonical identification
$\Det_{x_0}(M;\rho)=\mathbb C$, and via this identification 
$\mathcal T^{\e,\mathfrak o}_\comb(\rho)=(\tau_\comb^{\rho,\smash{\e},\mathfrak o})^{-1}$.

The function $\mathcal T_\comb^{\e,\mathfrak o}$ contains relevant
topological information, even in the simplest possible examples. In
\fullref{SS:maptor} we compute this function for a mapping torus in
terms of the Lefschetz zeta function of the gluing diffeomorphism.

If $M$ is obtained by surgery on a framed knot and $\dim V=1$, the function
$\mathcal T_\comb^{\e,\mathfrak o}$ coincides with the Alexander polynomial
of the knot; see Turaev \cite{Tu02}.

If $\e_1,\e_2\in\Eul_{x_0}(M;\Z)$ are two Euler structures with $P(\e_1)=P(\e_2)$ and 
$\mathfrak o_1$, $\mathfrak o_2$ are two cohomology orientations then the two functions 
$\mathcal T^{\e_1,\mathfrak o_1}_\comb$ and $\mathcal T^{\e_2,\mathfrak o_2}_\comb$ 
differ by multiplication with a root of unity; see \eqref{E:dep}. 
Therefore, up to multiplication by a root of unity,
$\mathcal T_\comb^{\e,\mathfrak o}$ depends only on $\e^*=P(\e)$.   
The Bismut--Zhang theorem implies that the absolute value of 
$\mathcal T_\comb^{\e,\mathfrak o}$ actually coincides with the 
modified Ray--Singer torsion associated with $\e^*$, \ie the quantity
$$
\mathcal T_\an^{\e^*}(\rho)=
T_\an(\nabla_\rho,g,\mu)\cdot e^{-\int_M\omega(\nabla_\rho,\mu)\wedge\alpha}.
$$
Here $\mathbb F_\rho=(F_\rho,\nabla_\rho)$ denotes the flat complex vector 
bundle associated to the
representation $\rho\in\Rep^M_o(\Gamma;V)$, $g$ is a Riemannian metric on
$M$, $\mu$ is a Hermitian structure on $F_\rho$, $T_\an(\nabla_\rho,g,\mu)$
denotes the Ray--Singer torsion, $\omega(\nabla_\rho,\mu)$ is the
Kamber--Tondeur form (\fullref{SS:KT}) and $\alpha$ is such 
that the pair $(g,\alpha)$ represents the co-Euler structure $\e^*$.
Recall that for acyclic $\rho$ we have a canonical identification 
$\Det_{x_0}(M;\rho)=\mathbb C$ and via this identification
$\mathcal T^{\e^*}_\an(\rho)=\|1\|_{\an}^{\rho,\e^*}$.

We use this fact to conclude (under some additional hypotheses) that  
the Ray--Singer torsion, when viewed as a function on the 
analytic set of flat connections, is the absolute value of a meromorphic 
function. The precise statement is contained in \fullref{C:meroa}.
This partially answers the question formulated at the beginning of this
introduction. This was also known to W\,M\"uller at least restricting to unimodular flat connections. A more analytic treatment is 
announced in \cite{BH05,BH05a,BH06}; see also \cite{BK05a,BK05,BK05b,BK06}.

The argument of the complex number $\mathcal T_\comb^{\e,\mathfrak o}(\rho)$  
is an apparently unexplored invariant. If $U$ is an open neighborhood of
$[1,2]\times\{0\}\subseteq\C$ and $\rho\co U\to\Rep_0^M(\Gamma;V)$ is a
holomorphic path in the nonsingular part of $\Rep_0^M(\Gamma;V)$
joining $\rho_1=\rho(1)$ with $\rho_2=\rho(2)$, then the difference 
$$
\arg\mathcal T_\comb^{\e,\mathfrak o}(\rho_2)-\arg\mathcal T_\comb^{\e,\mathfrak o}(\rho_1)
\in\R/2\pi\Z
$$
can be derived from $\mathcal T^{\e^*}_{\an}(\rho(z))$, $P(\e)=\e^*$; see
\eqref{E:ph} in \fullref{SS:inv}.
In general two representations $\rho_1$ and $\rho_2$ in the same 
connected component of the nonsingular part of $\Rep^M _0(\Gamma;V)$ 
cannot be joined by a holomorphic path. However, one can 
always find a finite collection of representations
$\rho^i\in\Rep^M_{\smash{0}}(\Gamma;V)$, $i=1,\dotsc,k$, with $\rho^1=\rho_1$ and $\rho^k=\rho_2$
so that holomorphic paths from $\rho^i$ to $\rho^{i+1}$ do exist, 
and implicitly benefit from such formula.

As another application we provide a derivation of Marcsik's theorem 
(unpublished \cite{M98}) from the Bismut--Zhang theorem and the computation of 
$\mathcal T_\comb^{\e,\mathfrak o}$ for mapping tori; see \fullref{SS:marcsik}.

{\bf Acknowledgements}\qua
We would like to thank V\,Turaev for pointing out a sign mistake
and some lacking references in a previous version of the paper.
We also would like to thank the referee for suggestions and requests of clarifications. 
We believe they have improved the quality of the exposition.
Part of this work was done while both authors enjoyed the
hospitality of the Max Planck Institute for Mathematics in Bonn.
A previous version was written while the second author enjoyed the
hospitality of the Ohio State University.
The second author was partially supported by the \emph{Fonds zur F\"orderung der
wissenschaftlichen Forschung} (Austrian Science Fund), project number {\tt P14195-MAT}.

\section[A few characteristic forms and the invariant R]{A few characteristic forms and the invariant $R$}\label{char}

\subsection{Euler, Chern--Simons, and the Mathai--Quillen form}\label{SS:MQECS}

In this section we will briefly recall basic properties of the Euler 
form, the Mathai--Quillen form and the Chern--Simons class. For more 
details we refer to \cite{MQ86} or \cite[section~III]{BZ92}.
Although these forms are defined for any real vector bundle equipped
with a connection and a parallel Hermitian metric we will only discuss them
for the tangent bundle of a Riemannian manifold and the Levi-Civita connection.

Let $M$ be smooth closed manifold of dimension $n$. Let $\pi\co TM\to M$ denote the
tangent bundle and $\Or_M$ the orientation bundle, a flat real line bundle over $M$. 
For a Riemannian metric $g$ denote its Euler form by 
\begin{gather*}
E(g)\in\Omega^n(M;\Or_M)\\ 
\Psi(g)\in\Omega^{n-1}(TM\setminus M;\pi^*\Or_M)\tag*{\hbox{and by}}
\end{gather*}
its Mathai--Quillen form. For two Riemannian metrics $g_1$ and $g_2$ we let 
$$
\cs(g_1,g_2)\in\Omega^{n-1}(M;\Or_M)/d(\Omega^{n-2}(M;\Or_M))
$$ 
denote their Chern--Simons class. The following relations hold:
\begin{eqnarray}
d\cs(g_1,g_2) &=& E(g_2)-E(g_1)
\label{E:csg:i}
\\
\cs(g_2,g_1) &=& -\cs(g_1,g_2)
\label{E:csg:ii}
\\
\cs(g_1,g_3) &=& \cs(g_1,g_2) + \cs(g_2,g_3)
\label{E:csg:iii}
\\
E(g)&=&(-1)^nE(g) 
\label{E:oddnE}
\\
\cs(g_1,g_2)&=&(-1)^n\cs(g_1,g_2)
\label{E:oddncs}
\end{eqnarray}
When $n$ is odd both, the Euler form and the Chern--Simons class, vanish.
Moreover, we have the equalities
\begin{eqnarray}
d\Psi(g) &=& \pi^*E(g)
\label{dpsie}
\\\label{pg1pg2cs}
\Psi(g_2)-\Psi(g_1) &\equiv& \pi^*\cs(g_1,g_2) 
\mod d\Omega^{n-2}(TM\setminus M;\pi^*\Or_M)
\\
\nu^*\Psi(g)&=&(-1)^n\Psi(g)
\label{E:nupsi}
\end{eqnarray}
where $\nu\co TM\to TM$ denotes the canonical involution, $\nu(x):=-x$.

For $M=\R^n$ equipped with the standard Riemannian metric $g_0$ we have
\begin{equation}\label{E:psig0}
\Psi(g_0)
=\frac{\Gamma(n/2)}{(2\pi)^{n/2}}\sum_{i=1}^n(-1)^i
\frac{\xi_i}{(\sum \xi_i^2)^{n/2}}
d\xi_1\wedge\cdots\wedge\widehat{d\xi_i}\wedge\cdots\wedge d\xi_n
\end{equation}
in standard coordinates $x_1,\dotsc,x_n,\xi_1,\dotsc,\xi_n$ on $TM$.
In general, the restriction of $\Psi$ to a fiber $T_xM$
coincides with the negative of the standard generator of the group
$H^{n-1}(T_xM\setminus0;\mathcal O_{T_xM})$.

Suppose $X$ is a vector field with isolated zero $x$, and let
$B_\epsilon(x)$ denote the ball of radius $\epsilon$ centered at $x$,
with respect to some chart. From \eqref{pg1pg2cs} and 
\eqref{E:psig0} one easily deduces
\begin{equation}\label{psiind}
\lim_{\epsilon\to0}\int_{\partial(M\setminus B_\epsilon(x))}
X^*\Psi(g)=\IND_X(x),
\end{equation}
where $\IND_X(x)$ denotes the Hopf index of $X$ at $x$.

\subsection{Euler and Chern--Simons class for vector fields}\label{SS:ECSX}

Let $M$ be a closed manifold and suppose $X$ is a vector field with isolated
set of zeros $\mathcal X\subseteq M$. For every $x\in\mathcal X$ we have a
Hopf index $\IND_X(x)$ of $X$ at $x$. Define a singular\footnote {Here, and
throughout the paper, singular chain stands for \emph{smooth\/} singular chain.}
zero-chain
$$
E(X):=\sum_{x\in\mathcal X}\IND_X(x)x\in C_0(M;\Z).
$$
Given two vector fields $X_1$ and $X_2$ we are going to define a singular 
one-chain, well-defined up to boundaries,
\begin{equation}\label{E:csx}
\cs(X_1,X_2)\in C_1(M;\Z)/\partial C_2(M;\Z)
\end{equation}
with the properties:
\begin{eqnarray}
\partial\cs(X_1,X_2) &=& E(X_2)-E(X_1)
\label{E:csx:i}
\\
\cs(X_2,X_1) &=& -\cs(X_1,X_2)
\label{E:csx:ii}
\\
\cs(X_1,X_3) &=& \cs(X_1,X_2) + \cs(X_2,X_3)
\label{E:csx:iii}
\\
E(-X)&=&(-1)^nE(X)
\label{E:oddnex}
\\
\cs(-X_1,-X_2)&=&(-1)^n\cs(X_1,X_2)
\label{E:oddncsx}
\end{eqnarray}
This should be compared with the equalities
\eqref{E:csg:i}--\eqref{E:oddncs} in \fullref{SS:MQECS}.

Equation \eqref{E:oddnex} follows from $\IND_{-X}(x)=(-1)^n\IND_X(x)$.
The construction of \eqref{E:csx} is accomplished by first reducing to
vector fields with nondegenerate zeros with the help of a small perturbation,
and then taking the zero set of a generic homotopy connecting $X_1$ with $X_2$.

More precisely, suppose we have two vector fields $X_1'$ and $X_2'$ 
with nonde\-ge\-ne\-rate zeros. 
Consider the vector bundle $p^*TM\to I\times M$, where $I$ is the interval $[1,2]$ and 
$p\co I\times M\to M$ denotes the natural projection. Choose a section
$\mathbb X$ of $p^*TM$ transversal to the zero section and which
restricts to $X_i'$ on $\{i\}\times M$, $i=1,2$. The zero set of $\mathbb X$
is a canonically oriented one-dimensional submanifold with boundary. Its
fundamental class, when pushed forward via $p$, gives rise to
$c(\mathbb X)\in C_1(M;\Z)/\partial C_2(M;\Z)$ satisfying \eqref{E:csx:i}.

Suppose $\mathbb X_1$ and $\mathbb X_2$ are two nondegenerate homotopies 
from $X_1'$ to $X_2'$. Then we have
$c(\mathbb X_1)=c(\mathbb X_2)\in C_1(M;\Z)/\partial C_2(M;\Z)$.
Indeed, consider the vector bundle
$q^*TM\to I\times I\times M$, where $q\co I\times I\times M\to M$ denotes the 
natural projection. Choose a section of $q^*TM$ which is transversal to the
zero section, restricts to $\mathbb X_i$ on $\{i\}\times I\times M$,
$i=1,2$, and restricts to $X_i'$ on $\{s\}\times\{i\}\times M$ for all $s\in I$
and all $i=1,2$. The zero set of such a section then gives rise to $\sigma$
satisfying $c(\X_2)-c(\X_1)=\partial\sigma$. So we get a well-defined 
$\cs(X_1',X_2'):=c(\X)\in C_1(M;\Z)/\partial C_2(M;\Z)$ satisfying
\eqref{E:csx:i}, provided $X_1'$ and $X_2'$ have nondegenerate zero set.

Now suppose $X$ is a vector field with isolated zeros. For every zero 
$x\in\mathcal X$ we choose a contractible open neighborhood
$B_x$ of $x$, assuming all $B_x$ are disjoint. Set
$B:=\bigcup_{x\in\mathcal X}B_x$. Choose a vector field with nondegenerate
zeros $X'$ that coincides with $X$ on $M\setminus B$. Let $\mathcal X'$
denote its zero set.  Then
$\IND_X(x)=\sum_{y\in\mathcal X'\cap B_x}\IND_{X'}(y)$ for every $x\in\mathcal X$.
So we can choose a one-chain $c(X,X')$ in $B$ which
satisfies $\partial c(X,X')=E(X')-E(X)$. Since $H_1(B;\Z)$ vanishes, the
one-chain $c(X,X')$ is well-defined up to a boundary.

Given two vector fields $X_1$ and $X_2$ with isolated zeros we choose
$B$, $X_1'$, $X_2'$, $c(X_1,X_1')$ and $c(X_2,X_2')$ as above and set 
$$
\cs(X_1,X_2):=c(X_1,X_1')+\cs(X_1',X_2')-c(X_2,X_2').
$$
Note that for two vector fields which coincide on $M\setminus B$ and 
have nondegenerate zeros only, one can construct a nondegenerate 
homotopy between them which is constant on $M\setminus B$. Together
with $H_1(B;\Z)=0$ this immediately implies that the one-chain
$\cs(X_1,X_2)\in C_1(M;\Z)/\partial C_2(M;\Z)$ is independent 
of these choices. Clearly \eqref{E:csx:i}--\eqref{E:csx:iii} and
\eqref{E:oddncsx} hold too.

\subsection{Kamber--Tondeur one-form}\label{SS:KT}

Let $E$ be a real or complex vector bundle over $M$.
For a connection $\nabla$ and a Hermitian structure $\mu$ on $E$ 
we define a real-valued one-form 
$\omega(\nabla,\mu):=-\frac12\tr_\mu(\nabla\mu)$ in $\Omega^1(M;\R)$.
More explicitly, for a tangent vector $X\in TM$
$$
\omega(\nabla,\mu)(X):=-\frac12\tr_\mu(\nabla_X\mu).
$$
This form is known as \emph{Kamber--Tondeur\/} form \cite{KT}. Our
convention differs from the one in Bismut and Zhang \cite{BZ92} by a factor $-\frac12$.
Here one interprets  $\mu$ as a section in the bundle $H$ of sesquilinear 
forms on $E$ equipped with the induced connection (still denoted by $\nabla$) 
and one denotes  by $\tr_\mu$ the bundle map 
from $H$ to the trivial rank one bundle defined by $\mu$. 
Note that if $\nabla$ is flat then $\omega(\nabla,\mu)$ will be 
a closed one-form. This and the following formulas are easily seen by 
noticing that $\omega(\nabla^{\det E},\mu^{\det E})=\omega(\nabla,\mu)$
where $\nabla^{\det E}$ and $\mu^{\det E}$ denote the induced connection
and Hermitian structure on the determinant line 
$\det E:=\Lambda^{\rank E}E$. The verification of this equality is 
straightforward. One calculates $\omega$ in coordinates and with 
respect to a $\mu$--orthonormal frame in the neighborhood of  
a given point $x\in M$. A simple comparison between the two sides 
of the equality above verifies the statement.

Recall that the space of connections on $E$ is an affine space over
$\Omega^1(M;E^*\otimes E)$. 
For two connections $\nabla^1$ and $\nabla^2$ we find
\begin{equation}\label{E:onn}
\omega(\nabla^2,\mu)-\omega(\nabla^1,\mu)=\Re(\tr(\nabla^2-\nabla^1))
\end{equation}
where the trace is with respect to the last two variables, and $\Re$
denotes the real part.

If $\mu_1$ and $\mu_2$ are Hermitian structures then there is
a unique section $A\in C^\infty(E^*\otimes E)$ such that
$\mu_2(e,f)=\mu_1(Ae,f)$. We set $V(\mu_1,\mu_2):=|\det A|$ 
which is a positive real-valued function on $M$. Its value at $x\in M$
is the volume, with respect to $(\mu_2)_x$, of a parallelepiped
obtained from an orthonormal base, with respect to $(\mu_1)_x$.
We have
\begin{equation}\label{E:omm}
\omega(\nabla,\mu_2)-\omega(\nabla,\mu_1)=-\frac12d\log V(\mu_1,\mu_2).
\end{equation}

Suppose $\nabla$ is flat, and let $U\subseteq M$ be a contractible open set.
Then we find a parallel Hermitian structure $\tilde\mu$ over $U$, and can use 
\eqref{E:omm} to obtain a local primitive of the closed one-form
$\omega(\nabla,\mu)|_U=d(-\frac12\log V(\tilde\mu,\mu))$.

\begin{remark}\label{R:KT}
Suppose $\nabla$ is flat, and let $x_0$ be a base point in $M$. Parallel
transport defines a right action of $\pi_1(M,x_0)$ on $E_{x_0}$, the fiber
over $x_0$. Composing with the inversion we obtain a representation
$\rho\co\pi_1(M,x_0)\to\GL(E_{x_0})$. Denote by $[\det\circ\rho]\co H_1(M;\Z)\to\C^*$
the homomorphism induced by the homomorphism
$\det\circ\rho\co\pi_1(M,x_0)\to\C^*$. Let $\log|[\det\circ\rho]|\in
H^1(M;\R)$ denote the cohomology class corresponding to the homomorphism
$\log|[\det\circ\rho]|\co H_1(M;\Z)\to\R$. This class is called the
\emph{Kamber--Tondeur class\/} and is represented by the Kamber--Tondeur form
$\omega(\nabla,\mu)$ for every Hermitian metric $\mu$.
Put another way, the cohomology class of $\omega(\nabla,\mu)$ captures 
the absolute value of the holonomy in $(\det E,\nabla^{\det E})$.
\end{remark}

The following observation will be useful in the proof of \fullref{C:meroa}.

\begin{proposition}\label{P:hol_omega}
Let $E$ be a complex vector bundle over $M$ equipped with a Hermitian
structure $\mu$. Then there exist mappings
$\nabla\mapsto\tilde\omega(\nabla,\mu)$ which are affine
over the map $\tr\co\Omega^1(M;E^*\otimes E)\to\Omega^1(M;\C)$ such that
$\Re(\tilde\omega(\nabla,\mu))=\omega(\nabla,\mu)$ and such that
$\tilde\omega(\nabla,\mu)$ is closed whenever $\nabla$ is flat. 
\end{proposition}

\begin{proof}
Choose a reference connection $\nabla'$ on $E$, and define
$$
\tilde\omega(\nabla,\mu):=\tr(\nabla-\nabla')+\omega(\nabla',\mu).
$$ 
This is certainly affine over the trace, and its real part coincides with
$\omega(\nabla,\mu)$ in view of \eqref{E:onn}. 
If $E$ admits flat connections we can choose $\nabla'$ to be flat.
Then $\tilde\omega(\nabla,\mu)$ will be closed provided $\nabla$ was flat.
\end{proof}

\subsection{Regularization of an integral}\label{SS:reg}

Let $M$ be a closed manifold of dimension $n$. 
Suppose $\alpha\in\Omega^{n-1}(M\setminus\mathcal X;\Or_M)$ is a form
with isolated set of singularities $\mathcal X\subseteq M$, and
suppose $d\alpha\in\Omega^n(M;\Or_M)$ is globally smooth.
Here $\Or_M$ denotes the orientation bundle of $M$, a flat real line bundle
over $M$. For $x\in\mathcal X$ define a real number
$$
\IND_\alpha(x):=\lim_{\epsilon\to0}\int_{\partial(M\setminus B_\epsilon(x))}\alpha
=\int_{\partial(M\setminus B_{\epsilon'}(x))}\alpha
+\int_{B_{\epsilon'}(x)}d\alpha
$$
where $B_\epsilon(x)$ denotes the ball of radius $\epsilon$ in a chart
centered at $x$ and $\epsilon'$ is small enough. Note that for a Riemannian metric $g$, 
and a vector field $X$ with isolated zero $x$, the form $X^*\Psi(g)$ has an isolated
singularity at $x$, $dX^*\Psi(g)=E(g)$ is globally smooth, see
\eqref{dpsie}, and in view of \eqref{psiind} we have
\begin{equation}\label{E:INDaINDX}
\IND_{X^*\Psi(g)}(x)=\IND_X(x).
\end{equation}

Let $\alpha$ be as above, and suppose
$\omega\in\Omega^1(M)$ is a closed real or complex valued one-form. 
Choose a function $f$ such that $f$ vanishes on $\mathcal X$ and such that 
$\omega':=\omega-df$ vanishes locally around $\mathcal X$. Define
$$
S(\omega,\alpha)
:=S(\omega,\alpha;f)
:=\int_{M\setminus\mathcal X}\omega'\wedge\alpha-\int_Mfd\alpha.
$$
The notation is justified by the following lemma.

\begin{lemma}\label{L:S}
The quantity $S(\omega,\alpha;f)$ does not depend on the choice of $f$.
It is linear in $\omega$ and $\alpha$. It satisfies:
\begin{enumerate}
\item\label{S:i}
$S(dh,\alpha)=\sum_{x\in\mathcal X}h(x)\IND_\alpha(x)-\int_Mhd\alpha$
for every smooth function $h$.
\item\label{S:ii}
$S(\omega,\beta)=\int_M\omega\wedge\beta$ 
for all $\beta\in\Omega^{n-1}(M;\Or_M)$. Particularly,
$S(\omega,d\gamma)=0$ for all $\gamma\in\Omega^{n-2}(M;\Or_M)$.
\end{enumerate}
\end{lemma}

\begin{proof}
Let us first show that $S(\omega,\alpha;f)$ does not depend on $f$.
So suppose $f_i$ are functions which vanish on $\mathcal X$ and 
are such that $\omega_i':=\omega-df_i$ vanishes locally around $\mathcal X$.
Then $\omega_2'-\omega_1'=-d(f_2-f_1)$ and $f_2-f_1$ vanish locally around
$\mathcal X$. Using Stokes' theorem for forms with compact support on
$M\setminus\mathcal X$ we get
$$
S(\omega,\alpha;f_2)-S(\omega,\alpha;f_1)
=-\int_{M\setminus\mathcal X}d(f_2-f_1)\wedge\alpha-\int_M(f_2-f_1)d\alpha=0.
$$

Now let us turn to \itemref{S:i}.
Choose a function $f$ so that
$dh-df$ vanishes locally around $\mathcal X$ and so that $f$ vanishes on
$\mathcal X$. For every $x\in\mathcal X$ choose an embedded ball
$B_\epsilon(x)$ around $x$ such that $dh-df$ vanishes on 
$B_\epsilon(x)$. Then $h-f$ is constant equal to $h(x)$ on $B_\epsilon(x)$, 
and Stokes' theorem implies 
$$
S(dh,\alpha)=\int_{M\setminus\mathcal X}(dh-df)\wedge\alpha-\int_Mfd\alpha
=\sum_{x\in\mathcal X}h(x)
\int_{\partial(M\setminus B_\epsilon(x))}\alpha-\int_Mhd\alpha.
$$
With $\epsilon\to0$ the statement follows. The remaining assertions are obvious.
\end{proof}

\subsection[The invariant R]{The invariant $R$}\label{SS:R}

Let $M$ be a closed manifold of dimension $n$. For a closed one-form
$\omega$ and two Riemannian metrics $g_1$ and $g_2$ we set
\begin{equation}\label{E:rgg}
R(\omega,g_1,g_2):=\int_M\omega\wedge\cs(g_1,g_2).
\end{equation}
Although $\cs(g_1,g_2)$ (\fullref{SS:MQECS}), is only defined up to 
$d\Omega^{n-2}(M;\Or_M)$ the quantity \eqref{E:rgg} is unambiguously defined 
according to Stokes' theorem.
If the dimension of $M$ is odd then $R(\omega,g_1,g_2)=0$; see~\eqref{E:oddncs}.

For a closed one-form $\omega$, and two vector fields $X_1$ and $X_2$ with isolated 
zeros we set 
\begin{equation}\label{E:rxxx}
R(\omega,X_1,X_2):=\int_{\cs(X_1,X_2)}\omega.
\end{equation}
Again, even though $\cs(X_1,X_2)$ is only defined up to a 
boundary (\fullref{SS:ECSX}) the quantity \eqref{E:rxxx} is well-defined in view of Stokes' theorem.

For a Riemannian metric $g$ and a vector field $X$ with isolated zeros 
$\mathcal X\subseteq M$ we have $X^*\Psi(g)\in\Omega^{n-1}(M\setminus\mathcal X;\Or_M)$
and $dX^*\Psi(g)=E(g)$; see \fullref{SS:MQECS}, especially \eqref{dpsie}. 
So we may define 
\begin{equation}\label{E:rxx}
R(\omega,X,g):=S(\omega,X^*\Psi(g))
\end{equation}
for every closed one-form $\omega$; see \fullref{SS:reg}. The invariant $R(\omega,X,g)$ has first 
appeared in the work of Bismut and Zhang \cite{BZ92}.

\begin{proposition}\label{P:R}
On $M$, let $\omega$ be a closed one-form, $h$ a smooth function,
$g$ and $g_i$ Riemannian metrics, and let $X$ and $X_i$ be vector fields with 
isolated zeros. Then:
\begin{enumerate}
\item\label{P:R:o}
$R(\omega,g_1,g_2)$, $R(\omega,X_1,X_2)$ and $R(\omega,X,g)$ are linear in $\omega$.
\item\label{P:R:oo}
$R(\omega,g_2,g_1)=-R(\omega,g_1,g_2)$.
\item\label{P:R:ooo}
$R(\omega,X_2,X_1)=-R(\omega,X_1,X_2)$.
\item\label{P:R:i}
$R(dh,X,g)=h(E(X))-\int_MhE(g)$.
\item\label{P:R:ii}
$R(dh,g_1,g_2)=\int_MhE(g_1)-\int_MhE(g_2)$.
\item\label{P:R:iii}
$R(dh,X_1,X_2)=h(E(X_2))-h(E(X_1))$.
\item\label{P:R:iv}
$R(\omega,X_1,X_3)=R(\omega,X_1,X_2)+R(\omega,X_2,X_3)$.
\item\label{P:R:v}
$R(\omega,g_1,g_3)=R(\omega,g_1,g_2)+R(\omega,g_2,g_3)$.
\item\label{P:R:vi}
$R(\omega,X,g_2)-R(\omega,X,g_1)=R(\omega,g_1,g_2)$.
\item\label{P:R:vii}
$R(\omega,X_2,g)-R(\omega,X_1,g)=R(\omega,X_1,X_2)$.
\end{enumerate}
Here $h(E(X))=\sum_{x\in\mathcal X}h(x)\IND_X(x)$ where $\mathcal
X$ denotes the zero set of $X$.
\end{proposition}

\begin{proof}
Statements \itemref{P:R:o}--\itemref{P:R:ooo} and 
\itemref{P:R:ii}--\itemref{P:R:v} follow immediately from Stokes'
theorem, \fullref{L:S}, equations \eqref{E:csg:i}--\eqref{E:csg:iii}
and equations \eqref{E:csx:i}--\eqref{E:csx:iii}.

Assertion \itemref{P:R:i} follows from \fullref{L:S}\itemref{S:i},
equation \eqref{E:INDaINDX}, and from $dX^*\Psi(g)=E(g)$; see \eqref{dpsie}.

Statement \itemref{P:R:vi} follows from \fullref{L:S}\itemref{S:ii} and 
$X^*\Psi(g_2)-X^*\Psi(g_1)=\cs(g_1,g_2)$; see~\eqref{pg1pg2cs}. 
Statements \itemref{P:R:i} and \itemref{P:R:vi} can be also found in \cite{BZ92}.

It remains to prove \itemref{P:R:vii}. 
In view of \itemref{P:R:i} and \itemref{P:R:iii} we may assume that $\omega$
vanishes in a neighborhood of the zero set 
$\mathcal X_1\cup\mathcal X_2$. It is then clear that we may also assume
that $X_1$ and $X_2$ have nondegenerate zeros only. Choose a nondegenerate homotopy
$\X$ from $X_1$ to $X_2$. Perturbing the homotopy slightly, cutting it
into several pieces and using \itemref{P:R:iv} we may assume that the zero set 
$\X^{-1}(0)\subseteq I\times M$ is actually contained in a closed simply connected
$I\times V$. Again, we may assume that
$\omega$ vanishes on $V$. Then obviously $R(\omega,X_1,X_2)=0$. Moreover, in
this situation Stokes' theorem implies
\begin{align*}
R(\omega,X_2,g)-R(\omega,X_1,g)
&=\int_{M\setminus V}\omega\wedge X_2^*\Psi(g)
-\int_{M\setminus V}\omega\wedge X_1^*\Psi(g)
\\&=
\int_{I\times(M\setminus V)}d(p^*\omega\wedge\X^*\tilde p^*\Psi(g))
\\&=
-\int_{I\times(M\setminus V)}p^*(\omega\wedge E(g))=0.
\end{align*}
where $p\co I\times M\to M$ denotes the natural projection, and
$\tilde p\co p^*TM\to TM$ denotes the natural vector bundle homomorphism
over $p$. For the last calculation note that $d\X^*\tilde p^*\Psi(g)=p^*E(g)$
in view of \eqref{dpsie}, and that $\omega\wedge E(g)=0$ because of dimensional reasons.
\end{proof}

\begin{remark}
A similar definition of $R(\omega,X,g)$ works for any vector field $X$ 
with arbitrary zero set $\mathcal X:=\{x\in M\mid X(x)=0\}$ provided $\omega$ 
is exact when restricted to a sufficiently small neighborhood of $\mathcal X$. 
\end{remark}

\subsection{Extension of Chern--Simons theory}\label{SS:CS}

The material in this section will not be used in the rest of the paper,
but puts the observations in \fullref{P:R} in a more natural perspective.

Let $M$ be a closed manifold of dimension $n$. Equip $\Omega^k(M;\R)$
with the $C^\infty$--topology. The continuous linear functionals on
$\Omega^k(M;\R)$ are called $k$--currents. Let $\mathcal D_k(M)$ denote the space of all $k$--currents and $\delta\co\mathcal D_k(M)\to\mathcal
D_{k-1}(M)$ be given by 
$$
(\delta\varphi)(\alpha):=\varphi(d\alpha).
$$

Clearly $\delta^2=0$. By the deRham theorem for currents the chain complex 
$(\mathcal D_*(M),\delta)$ computes the homology of $M$ with real coefficients.

We have a morphism of chain complexes
$$
C_*(M;\R)\to\mathcal D_*(M),\quad
\sigma\mapsto\hat\sigma,\quad
\hat\sigma(\alpha):=\int_\sigma\alpha.
$$
Here $C_*(M;\R)$ denotes the space of singular chains with real
coefficients. Moreover, we have a morphism of chain complexes
$$
\Omega^{n-*}(M;\Or_M)\to\mathcal D_*(M),\quad
\beta\mapsto\hat\beta,\quad
\hat\beta(\alpha):=(-1)^{\frac12|\alpha|(|\alpha|+1)}\int_M\alpha\wedge\beta.
$$
Here $|\alpha|$ denotes the degree of $\alpha$.

Note that any zero-current represents a degree zero homology class, and 
$H_0(M;\R)$ identifies canonically to $\bigoplus_\alpha\R_\alpha$
where $\R_\alpha$ is a copy of $\R$ corresponding to a 
connected component of $M=\bigsqcup_\alpha M_\alpha$. One denotes by
$\overline{\chi}(M)$ the vector
$\{\chi(M_\alpha)\}\in\bigoplus_\alpha\R_\alpha$, and by 
$\mathcal E_0$ the affine subspace of $\mathcal D_0(M)$ consisting 
of zero-currents which represent $\overline{\chi}(M)$. Elements in
$\mathcal E_0$ are called \emph{Euler currents\/}.

Every vector field $X$ with isolated set of zeros $\mathcal X_\tau$ gives rise 
to a zero-chain $E(X)$; see \fullref{SS:ECSX}. Via the first morphism we
get a zero-current $\hat E(X)$. More explicitly,
$\hat E(X)(h)=\sum_{x\in\mathcal X}\IND_x(X)h(x)$
for all functions $h\in\Omega^0(M;\R)$. 
By Hopf's theorem this is an Euler current.

A Riemannian metric $g$ has an Euler form $E(g)\in\Omega^n(M;\Or_M)$;
see \fullref{SS:MQECS}. Via the second morphism we get 
a zero-current $\hat E(g)$ which, by Gauss--Bonnet, is an Euler current. 
More explicitly, $\hat E(g)(h)=\int_MhE(g)$ for all functions $h\in\Omega^0(M;\R)$.

Let $\mathcal Z^k(M;\R)\subseteq\Omega^k(M;\R)$ denote the space of closed 
$k$--forms on $M$ equipped with the $C^\infty$--topology. The continuous linear
functionals on $\mathcal Z^k(M;\R)$ are referred to as \emph{$k$--currents 
rel.\ boundary} and identify to $\mathcal D_k(M)/\delta(\mathcal D_{k+1}(M))$.
The two chain morphisms provide mappings
\begin{gather}\label{crbi}
C_k(M;\R)/\partial(C_{k+1}(M;\R))\to
\mathcal D_k(M)/\delta(\mathcal D_{k+1}(M))\\
\label{crbii}
\Omega^{n-k}(M;\Or_M)/d(\Omega^{n-k-1}(M;\Or_M))\to
\mathcal D_k(M)/\delta(\mathcal D_{k+1}(M)).
\end{gather}

For two vector fields $X_1$ and $X_2$ with isolated zeros, in \fullref{SS:ECSX} we have
constructed the one-chain $\cs(X_1,X_2)\in C_1(M;\Z)/\partial(C_2(M;\Z))$ which then gives rise to 
$\cs(X_1,X_2)\in C_1(M;\R)/\partial(C_2(M;\R))$ and via \eqref{crbi} we get a
one-current rel.\ boundary which we will denote by $\hat\cs(X_1,X_2)$. More
precisely, for all closed one-forms $\omega\in\mathcal Z^1(M;\R)$, we have
$\hat\cs(X_1,X_2)(\omega)=\int_{\cs(X_1,X_2)}\omega=R(\omega,X_1,X_2)$.

For two Riemannian metrics $g_1$ and $g_2$, we have the Chern--Simons form
$\cs(g_1,g_2)$ in $\Omega^{n-1}(M;\Or_M)/d(\Omega^{n-2}(M;\Or_M))$; see
\fullref{SS:MQECS}.
Via \eqref{crbii} we get a one-current rel.\ boundary which we denote by
$\hat\cs(g_1,g_2)$. More precisely, for closed one-forms $\omega\in\mathcal Z^1(M;\R)$,  
we have $\hat\cs(g_1,g_2)(\omega)=-\int_M\omega\wedge\cs(g_1,g_2)=-R(\omega,g_1,g_2)$.

Suppose $X$ is a vector field with isolated zeros and $g$ is a Riemannian
metric. We define a one-current rel.\ boundary by 
$\hat\cs(X,g)(\omega):=-R(\omega,X,g)$ for all closed one-forms
$\omega\in\mathcal Z^1(M;\R)$; see \fullref{SS:R}. Moreover
set $\hat\cs(g,X):=-\hat\cs(X,g)$.

\fullref{P:R} can now be reformulated.

\begin{proposition}
Let any of the symbols $x,y,z$ denote either a Riemannian metric or a 
vector field with isolated zeros. Then one has:
\begin{enumerate}
\item 
$\delta\hat\cs(x,y)=\hat E(y)-\hat E(x)$.
\item 
$\hat\cs(y,x)=-\hat\cs(x,y)$.
\item 
$\hat\cs(x,z)=\hat\cs(x,y)+\hat\cs(y,z)$.
\end{enumerate}
\end{proposition}

\subsection{Smooth triangulations}\label{SS:triang}

Smooth triangulations provide a remarkable source of vector fields with 
isolated zeros. To any smooth triangulation $\tau$ of a smooth manifold $M$ one 
can associate a smooth vector field $X_\tau$ called \emph{Euler vector field\/},
with the following properties:

\begin{enumerate}
\item[P1]
The zeros of $X_\tau$ are all nondegenerate and are exactly the 
barycenters $x_\sigma$ of the simplices $\sigma$ of $\tau$.
\item[P2]
The piecewise differentiable function $f_\tau\co M\to\R$ defined by 
$f_\tau(x_\sigma)=\dim(\sigma)$, $x_\sigma$ the barycenter of the 
simplex $\sigma$ and extended 
by linearity on each simplex of the barycentric subdivision of $\tau$, is a
Lyapunov function for $-X_\tau$, \ie strictly decreasing on nonconstant
trajectories of $-X_\tau$.
\item[P3] 
For each zero $x_\sigma$ the unstable set with respect to $-X_\tau$ 
agrees in some small neighborhood of  $x_\sigma$
with the open simplex $\sigma$, consequently the zeros are hyperbolic. 
The Morse index of $-X_\tau$ at $x_\sigma$ equals $\dim(\sigma)$ and the 
(Hopf) index of $X_\tau$ at $x_\sigma$ equals $(-1)^{\dim(\sigma)}$. 
\end{enumerate}

The space of vector fields satisfying P1--P3 is a contractible space.
The convex combination provides a contraction.
It was pointed out to us that a somewhat similar vector field has been considered
by Steenrod \cite[Section 39.7]{S99}.

To construct such a vector field we first construct a piecewise smooth 
(actually Lipschitz) vector field $X$ which satisfies P1, P2 and a 
stronger form of P3. Precisely the unstable sets of the rest points 
identify to the open simplices of the triangulation.\footnote{It is not 
always possible to find a smooth vector field with such properties.}

We begin with a standard simplex $\Delta_n$ of vectors 
$(t_0,\dotsc,t_n)\in\R^{n+1}$ which satisfy $0\leq t_i\leq 1$ and $\sum t_i=1$. 
\begin{enumerate}
\item 
Let $E_n$ denote the Euler vector field of the corresponding affine space,
$\sum t_i=1$, centered at the barycenter $O=(1/(n+1),\dotsc,1/(n+1))$ 
and restricted to $\Delta_n$. 
\item
Let $e\co\Delta_n\to [0,1]$ denote the function which is $1$ on the barycenter
$O$, zero on the boundary and extended linearly on any ray through the
barycenter $O$.
\item 
Let $r\co\Delta_n\setminus\{O\}\to\partial\Delta_n$ denote the radial 
retraction to the boundary.
\end{enumerate}

Set $X'_n:=e\cdot E_n$, which is a vector field on $\Delta_n$.

By induction we will construct a canonical vector field $X_n$ on $\Delta_n$ 
which at any point $x\in\Delta_n$ is tangent to the open face the point 
belongs to and vanishes only at the barycenter of each face. We proceed as follows:

Suppose we have constructed such canonical vector fields on all $\Delta(k)$,
$k\leq n-1$. Using the canonical vector fields $X_{n-1}$ we define the vector
field $X_n$ on the boundary $\partial\Delta_n$ and extend it to the vector
field $X''_n$ by taking at each point $x\in\Delta_n$ the vector parallel to
$X_n(r(x))$ multiplied by the function $(1-e)$ and at $O$ the vector zero. 
Clearly, such a vector field vanishes on the radii $\overline{OP}$ (P the
barycenter of any face). We finally put
$$
X_n:=X'_n+X''_n.
$$
The vector field $X_n$ is continuous and piecewise differential 
(actually Lipschitz) and has a well-defined continuous flow.

Putting together the vector fields $X_n$ on all simplices we provide
a piecewise differential (and Lipschitz) vector field $X$ on any simplicial
complex or polyhedron and in particular on any smoothly
triangulated manifold. The vector field $X$ has a flow and $f_\tau$ is a
Lyapunov function for $-X$. The vector field $X$ is not necessary smooth but
by a small (Lipschitz) perturbation we can approximate it by a smooth vector
field which satisfies P1--P3. Any of the resulting vector fields is
referred to as the Euler vector field of a smooth triangulation $\tau$ and
denoted by $X_\tau$.

\begin{remark}\label{R:vect}
A similar construction can be done for smooth cell structures. Recall that
a smooth cell structure on a smooth manifold $M$ is a CW-complex structure
provided by cells whose characteristic maps 
$\chi\co D^k\to M$ are smooth embeddings (one-to-one and of maximal rank at 
any point including the corners) of compact $k$--dimensional convex sets 
$D$, with the property that the image of the boundary is a union 
of open cells of lower dimension.\footnote{More general cells can be 
considered but for the need of this paper (\fullref{SS:maptor}) this suffices.} 
We are led to a well-defined homotopy class of smooth vector fields  
whose rest points are the barycenters of the cells, are nondegenerate
and of the same index as the dimension of the cell.
\end{remark}

\section{Euler and co-Euler structures}\label{S:eul}

\subsection{Euler structures}\label{SS:eul}

Euler structures have been introduced by Turaev \cite{Tu90} for 
manifolds with vanishing Euler--Poincar\'e characteristic. 
In the presentation below we remove the vanishing
hypothesis at the expense of a base point.

Let $(M,x_0)$ be a base pointed closed connected manifold of dimension $n$.
Let $X$ be a vector field and let $\mathcal X$ denote its zero set.
Suppose the zeros of $X$ are isolated and recall the singular zero-chain
$E(X)=\sum_{x\in\mathcal X}\IND_X(x)x$ from \fullref{SS:ECSX}.
An \emph{Euler chain\/} for $X$ is a singular one-chain $c\in C_1(M;\Z)$ so that
$$
\partial c=E(X)-\chi(M)x_0.
$$
Since $\sum_{x\in\mathcal X}\IND_X(x)=\chi(M)$ every vector field with 
isolated zeros admits Euler chains.

Consider pairs $(X,c)$ where $X$ is a vector field with isolated zeros,
and $c$ is an Euler chain for $X$. We call two such pairs $(X_1,c_1)$ and
$(X_2,c_2)$ equivalent if
$$
c_2=c_1+\cs(X_1,X_2)\in C_1(M;\Z)/\partial(C_2(M;\Z)).
$$
For the definition of $\cs(X_1,X_2)$ see \fullref{SS:ECSX}. We will write
$\Eul_{x_0}(M;\Z)$ for the set of equivalence classes as above and 
$[X,c]\in\Eul_{x_0}(M;\Z)$ for the element represented by the pair $(X,c)$. 
Elements of $\Eul_{x_0}(M;\Z)$ are called \emph{Euler structures of $M$ based at 
$x_0$.} There is an obvious $H_1(M;\Z)$ action on $\Eul_{x_0}(M;\Z)$ defined by
$$
[X,c]+[\sigma]:=[X,c+\sigma],
$$
for $[\sigma]\in H_1(M;\Z)$ and $[X,c]\in\Eul_{x_0}(M;\Z)$. Evidently this action 
is free and transitive. In this sense $\Eul_{x_0}(M;\Z)$ is an affine version
of $H_1(M;\Z)$.

Considering Euler chains with real coefficients 
one obtains, in exactly the same way, an affine version of $H_1(M;\R)$ 
which we will denote by $\Eul_{x_0}(M;\R)$. Moreover, there is a natural map
$\Eul_{x_0}(M;\Z)\to\Eul_{x_0}(M;\R)$ which is affine over the homomorphism
$H_1(M;\Z)\to H_1(M;\R)$.

\begin{remark}\label{reeb_surgery}
There is an alternative way of understanding the $H_1(M;\Z)$ action on
$\Eul_{x_0}(M;\Z)$.
Suppose $n>2$ and represent $[\sigma]\in H_1(M;\Z)$ by a simple closed curve
$\sigma$. Choose a tubular neighborhood $N$ of $\sigma(S^1) $ considered as vector
bundle $N\to S^1$. Choose a fiber metric and a linear connection on $N$. 
Choose a representative of $[X,c]\in\Eul_{x_0}(M;\Z)$ such that 
$X|_N=\frac\partial{\partial\theta}$, the horizontal lift of the canonical 
vector field on $S^1$. Choose a function $\lambda\co[0,\infty)\to[-1,1]$,
which satisfies $\lambda(r)=-1$ for $r\leq\frac13$ and $\lambda(r)=1$ for
$r\geq\frac23$. Finally, choose a function $\mu\co[0,\infty)\to\R$ satisfying
$\mu(r)=r$ for $r\leq\frac13$, $\mu(r)=0$ for $r\geq\frac23$ and $\mu(r)>0$
for all $r\in(\frac13,\frac23)$. Now construct a new smooth vector field 
$\tilde X$ on $M$ by setting
$$
\tilde X:=
\begin{cases}
X & \text{\rm on $M\setminus N$}
\\
\lambda(r)\frac\partial{\partial\theta}
+\mu(r)\frac\partial{\partial r} 
& \text{\rm on $N$}
\end{cases}
$$
where $r\co N\to[0,\infty)$ denotes the radius function determined by the
fiber metric on $N$ and $-r\frac\partial{\partial r}$ is the Euler vector
field of $N$. This construction is known as Reeb surgery \cite{N03}.
If the zeros of $X$ are all nondegenerate the homotopy 
$X_t:=(1-t)X+t\tilde X$ is a nondegenerate homotopy from $X_0=X$ to
$X_1=\tilde X$ from which one easily deduces that
$[\tilde X,c]=[X,c]+[\sigma]$.
Particularly, all the choices that entered the Reeb surgery do not affect 
the out-coming Euler structure $[\tilde X,c]$.
\end{remark}

Let us consider a change of base point. Let $x_0,x_1\in M$ and choose a path
$\sigma$ from $x_0$ to $x_1$. Define
\begin{equation}\label{chbp}
\Eul_{x_0}(M;\Z)\to\Eul_{x_1}(M;\Z),\quad
[X,c]\mapsto[X,c-\chi(M)\sigma].
\end{equation}
This is an $H_1(M;\Z)$ equivariant bijection but depends on the homology
class of $\sigma$ rel.\ $\{x_0,x_1\}$. However, if $\chi(M)=0$, then 
$\Eul_{x_0}(M;\Z)$ does not depend on the base point \cite{Tu90},
and we will write $\Eul(M;\mathbb Z)$ in this case.

\begin{remark}\label{R:caneul}
The assignment $\nu([X,c]):=[-X,(-1)^nc]$ defines an involution $\nu$
on $\Eul_{x_0}(M;\mathbb R)$; see \eqref{E:oddnex} and \eqref{E:oddncsx}. 
It is affine over the homomorphism $(-1)^n$ from $H_1(M;\mathbb R)$ to 
$H_1(M;\mathbb R)$. Hence $\nu=\id$ for even $n$. If $n$ is odd
then $\nu$ has a unique fixed point which we denote by
$\e_\can\in\Eul(M;\mathbb R)$ and referred to as the \emph{canonical Euler
structure.} It can be represented as $\e_\can=[X,\frac12\cs(-X,X)]$
and this does not depend on $X$. Hence we obtain a canonical identification
$H_1(M;\mathbb R)=\Eul(M;\mathbb R)$, $\sigma\mapsto\e_\can+\sigma$,
provided $n$ is odd.
This involution makes sense on $\Eul_{x_0}(M;\mathbb Z)$ too,
but neither must it be trivial for even $n$, nor must it possess a (unique)
fixed point if $n$ is odd.
Finally, note that the natural mapping $\Eul_{x_0}(M;\Z)\to\Eul_{x_0}(M;\R)$
obviously intertwines the involutions on $\Eul_{x_0}(M;\Z)$ and $\Eul_{x_0}(M;\R)$.
\end{remark}

\begin{lemma}\label{vector-field}
Suppose $X$ is a smooth vector field with isolated zeros,
$\e\in\Eul_{x_0}(M;\Z)$ and assume $X$ has at least one zero with
Hopf index $\pm1$.
Then one can choose a collection of smooth paths $\pi_x\co[0,1]\to M$
with $\pi_x(0)=x_0$ and $\pi_x(1)=x$, $x\in\mathcal X$ so that for 
$c:=\sum_{x\in\mathcal X}\IND_X(x)\pi_x$ one has $\e=[X,c]$. 
Moreover, if $\dim M>2$ one can choose these paths to be smoothly embedded
and so that they only intersect at $x_0$.
\end{lemma}

\begin{proof}
Let $\pi_x'\co[0,1]\to M$ be any collection of smooth paths with
$\pi_x'(0)=x_0$ and $\pi'_x(1)=x$, $x\in\mathcal X$. Then
$c':=\sum_{x\in\mathcal X}\IND_X(x)\pi_x'$ is an Euler chain for
$X$, and we get an Euler structure $\e':=[X,c']\in\Eul_{x_0}(M;\mathbb Z)$.
Since the $H_1(M;\mathbb Z)$ action on $\Eul_{x_0}(M;\mathbb Z)$ is
transitive and 
the Hur\'ewicz homomorphism $\pi_1(M,x_0)\to H_1(M;\mathbb Z)$ is
onto, we find a closed smooth path $\sigma\co[0,1]\to M$,
$\sigma(0)=\sigma(1)=x_0$ satisfying $\e=\e'+[\sigma]$.
Choose a zero $y\in\mathcal X$ with $\IND_X(y)=\pm1$. Define $\pi_y$ as the concatenation
$\sigma^{\IND_X(y)}*\pi_y'$, and set $\pi_x:=\pi_x'$ for $y\neq x\in\mathcal X$.
Then $c:=\sum_{x\in\mathcal X}\IND_X(x)\pi_x$ is an Euler chain
for $X$ and $[X,c]=\e'+\IND_X(y)^2[\sigma]=\e$.

If $\dim(M)>2$ it is possible to approximate the paths $\pi_x$
by smoothly embedded paths $\tilde\pi_x$ which intersect only at $x_0$,
have the same endpoints, and such that the one-chain $\tilde\pi_x-\pi_x\in C_1(M;\mathbb
Z)$ is a boundary. Then $\tilde c:=\sum_{x\in\mathcal
X}\IND_X(x)\tilde\pi_x$ is an Euler chain for $X$ and $[X,\tilde c]=[X,c]=\e$. 
\end{proof}

Let $\vf(M,x_0)$ denote the space of vector fields which vanish at $x_0$ and
are nonzero elsewhere. We equip this space with the $C^\infty$--topology.
Let $\pi_0(\vf(M,x_0))$ denote the set of homotopy classes of such vector
fields. If $X\in\vf(M,x_0)$ we will write $[X]$ for the corresponding class
in $\pi_0(\vf(M,x_0))$. The following proposition is due to Turaev in the case
$\chi(M)=0$.

\begin{proposition}\label{euldb}
Suppose $n>2$. Then there exists a natural bijection
\begin{equation}\label{pi0vseul}
\pi_0(\vf(M,x_0))=\Eul_{x_0}(M;\Z),\qquad
[X]\mapsto[X,0].
\end{equation}
\end{proposition}

\begin{proof}
To see that \eqref{pi0vseul} is well-defined note that for
$[X_1]=[X_2]\in\pi_0(\vf(M,x_0))$ we find a representative of
$\cs(X_1,X_2)$ which is supported in a contractible open neighborhood
$U$ of $x_0$. Since $H_1(U;\mathbb Z)=0$ this representative is a 
boundary, hence $[X_1,0]=[X_2,0]\in\Eul_{x_0}(M;\mathbb Z)$.

To see that \eqref{pi0vseul} is onto let
$\e\in\Eul_{x_0}(M;\mathbb Z)$ be an Euler structure. Choose a vector field
$X$ on $M$ with isolated zero set $\mathcal X$ such that it has at
least one zero with Hopf index $\pm1$. In view of \fullref{vector-field}
there exists an Euler chain $c$ for $X$ with $[X,c]=\e$ and 
such that $c$ is supported in a compact contractible set containing $x_0$.
Hence we find a smoothly embedded
open disk $D\subseteq M$ so that $x_0\in D$, $\mathcal X\subseteq D$ 
and so that $c$ is supported in $D$.
Choose a vector field $X'$ which coincides with $X$ on a neighborhood
of $M\setminus D$ and
vanishes just at $x_0$. Then there exists a representative of 
$\cs(X,X')$ which is supported in $D$. Since $H_1(D;\Z)=0$ the difference
of this representative and $c$ must be a boundary and we obtain
$[X',0]=[X,c]=\e$. Therefore \eqref{pi0vseul} is onto.

To prove injectivity of \eqref{pi0vseul} let $X_1,X_2\in\vf(M,x_0)$ and
suppose $[X_1,0]=[X_2,0]\in\Eul_{x_0}(M;\mathbb Z)$, \ie
$\cs(X_1,X_2)=0\in H_1(M;\Z)$. Let $\bar D\subseteq M$ denote an embedded
closed disk (of full dimension) centered at $x_0$, and let $D$ denote its interior. 
Consider the vector bundle $p^*TM\to I\times M$
and consider the two vector fields as a nowhere vanishing section of $p^*TM$
defined over the set $\partial I\times\dot M$, where $\dot M:=M\setminus D$.
We would like to extend it to a nowhere vanishing section over
$I\times\dot M$. The first obstruction we meet is an element in
\begin{eqnarray*}
H^n(I\times\dot M,\partial I\times\dot M;\{\pi_{n-1}\})
&=&
H_1(I\times\dot M,I\times\partial D;\Z)
\\
&=&
H_1(M,\bar D;\Z)
\\
&=&
H_1(M;\Z)
\end{eqnarray*}
which corresponds to $\cs(X_1,X_2)=0$. Here $\{\pi_{n-1}\}$ denotes the system
of local coefficients determined by the sphere bundle of $p^*TM$ with 
$\pi_{n-1}= \pi_{n-1}(S^{n-1})$. Since this obstruction vanishes by
hypothesis the next obstruction is defined and is an element in
\begin{eqnarray*}
H^{n+1}(I\times\dot M,\partial I\times\dot M;\{\pi_n\})
&=&
H_0(I\times\dot M,I\times\partial D;\pi_n(S^{n-1}))
\\
&=&
H_0(M,\bar D;\pi_n(S^{n-1}))
\\
&=&
0.
\end{eqnarray*}
So obstruction theory as in Whitehead \cite{W78}
tells us that we find a nowhere vanishing section of $p^*TM$ defined over
$I\times\dot M$, which restricts to $X_i$ on $\{i\}\times\dot M$, $i=1,2$.
Such a section can easily be extended to a globally defined section of
$p^*TM\to I\times M$, which restricts to $X_i$ on $\{i\}\times M$, $i=1,2$,
and whose zero set is precisely $I\times\{x_0\}$.\footnote{For example
trivialize the bundle $p^*TM$ over $I\times\bar D$ and extend the section on
$\partial(I\times\bar D)$ to $I\times\bar D$ by scaling radially centered at 
$(\frac12,x_0)\in I\times \bar D$.} The resulting globally defined section can be
considered as homotopy from $X_1$ to $X_2$ showing
$[X_1]=[X_2]\in\pi_0(\vf(X,x_0))$. Hence \eqref{pi0vseul} is injective.
\end{proof}

\begin{remark}
Provided $n>2$, Reeb surgery defines an $H_1(M;\Z)$ action on the set $\pi_0(\vf(M,x_0))$
which via \eqref{pi0vseul} corresponds to the $H_1(M;\Z)$ action on
$\Eul_{x_0}(M;\Z)$; cf \fullref{reeb_surgery}.
\end{remark}

Let $\vf_0(M)$ denote the space of nowhere vanishing vector fields on $M$ 
equipped with the $C^\infty$--topology. Let $\pi_0(\vf_0(M))$ denote the
set of its connected components. For the next statement see also Turaev \cite{Tu90}.

\begin{proposition}\label{eulnv}
Suppose $\chi(M)=0$ and $n>2$. Then we have a surjection
\begin{equation}\label{pi00vseul}
\pi_0(\vf_0(M))\to\Eul_{x_0}(M;\Z),\qquad [X]\mapsto[X,0].
\end{equation}
\end{proposition}

\begin{proof}
The assignment \eqref{pi00vseul} is well-defined since for
$[X_1]=[X_2]\in\pi_0(\vf_0(M))$ the Chern--Simons class $\cs(X_1,X_2)$
can be represented by $0$.

To prove surjectivity let $\e\in\Eul_{x_0}(M;\mathbb Z)$ be an Euler structure. 
As explained in the proof of \fullref{euldb} we find a vector field
$X$ and an Euler chain $c$ for $X$ with $\e=[X,c]$ and such that there
exists a smoothly embedded open disk $D$ with $x_0\in D$, $\mathcal
X\subseteq D$ and so that $c$ is supported in $D$.
Wlog we may assume that $\bar D$ is a smoothly embedded closed disk in
$M$. Since $\chi(M)=0$ the degree of the map
$X\co\partial\bar D\to T\bar D\setminus\bar D$ vanishes. 
Hence we find a nowhere vanishing $X'$ which coincides with $X$ on 
a neighborhood of $M\setminus D$.
Hence $\cs(X,X')$ admits a representative supported in $D$. Since 
$H_1(D;\mathbb Z)=0$ the difference
of this representative and $c$ must be a boundary and we obtain
$[X',0]=[X,c]=\e$. Therefore \eqref{pi00vseul} is onto. 
\end{proof}

\subsection{Co-Euler structures}\label{SS:coeul}

We will now describe another approach to Euler structures which is in some 
sense Poincar\'e dual to the one in \fullref{SS:eul}. We still consider a closed
connected base pointed manifold $(M,x_0)$ of dimension $n$. Consider pairs
$(g,\alpha)$ where $g$ is a Riemannian metric on $M$ and 
$\alpha\in\Omega^{n-1}(M\setminus x_0;\mathcal O_M)$ with $d\alpha$ equal to the
Euler class $E(g) \in\Omega^n(M;\mathcal O_M)$ of $g$; see
\fullref{SS:MQECS}. Since $H^n(M\setminus x_0;\mathcal O_M)=0$ every Riemannian
metric admits such $\alpha$.
We call two pairs $(g_1,\alpha_1)$ and $(g_2,\alpha_2)$ equivalent if
$$
\cs(g_1,g_2)=\alpha_2-\alpha_1
\in\Omega^{n-1}(M\setminus x_0;\mathcal O_M)/
d\Omega^{n-2}(M\setminus x_0;\mathcal O_M).
$$
Write $\Eul_{x_0}^*(M;\R)$ for the set of equivalence classes and
$[g,\alpha]$ for the equivalence class represented by the pair $(g,\alpha)$.
Elements of $\Eul^*_{\smash{x_0}}(M;\R)$ are called \emph{co-Euler structures based at $x_0$\/}.
There is a natural $H^{n-1}(M;\mathcal O_M)$ action on $\Eul^*_{\smash{x_0}}(M;\R)$ given for $[\beta]\in H^{n-1}(M;\mathcal O_M)$ by
$$
[g,\alpha]+[\beta]:=[g,\alpha-\beta].
$$
Since 
$H^{n-1}(M;\mathcal O_M)=H^{n-1}(M\setminus x_0;\mathcal O_M)$ this action is
obviously free and transitive. In this sense $\Eul^*_{\smash{x_0}}(M;\R)$ is an affine
version of $H^{n-1}(M;\mathcal O_M)$.

\begin{remark}\label{metric}
Given a Riemannian metric $g$ and a class $\e^*\in\Eul^*_{\smash{x_0}}(M;\R)$ one can 
find $\alpha\in\Omega^{n-1}(M\setminus x_0;\mathcal O_M)$ with 
$d\alpha=E(g)$ so that $\e^*=[g,\alpha]$. Moreover, if $\chi(M)=0$ 
one can choose $\alpha\in\Omega^{n-1}(M;\mathcal O_M)$.
Hence in the case of vanishing Euler characteristic the set of
co-Euler structures is defined without reference to a base point
and we will denote it by $\Eul^*(M;\mathbb R)$.
\end{remark}

\begin{remark}\label{R:cancoeul}
Suppose $n$ is odd. Then $\chi(M)=0$, $E(g)=0$ and $\cs(g_1,g_2)=0$ for all
Riemannian metrics $g$, $g_1$ and $g_2$; see \eqref{E:oddnE} and
\eqref{E:oddncs}. Hence $[g,0]$ represents
a co-Euler structure independent of the Riemannian metric $g$. We refer
to this co-Euler structure as the \emph{canonical co-Euler structure\/} and will
denote it as $\e^*_\can\in\Eul^*(M;\mathbb R)$. It provides
a canonical identification $H^{n-1}(M;\mathcal O_M)=\Eul^*(M;\mathbb R)$,
$\beta\mapsto\e^*_\can+\beta$. This too can be understood in terms of the
involution $[g,\alpha]\mapsto[g,(-1)^n\alpha]$ on $\Eul^*_{\smash{x_0}}(M;\R)$; cf \fullref{R:caneul}.
\end{remark}

Suppose $[X,c]\in\Eul_{x_0}(M;\R)$, $[g,\alpha]\in\Eul^*_{\smash{x_0}}(M;\R)$
and let $\omega\in\Omega^1(M;\mathbb R)$ be a closed one-form.
Then the real number
\begin{equation}\label{rhoinv}
R(\omega,X,g)-S(\omega,\alpha)-\int_c\omega
\end{equation}
does only depend on the Euler structure $[X,c]\in\Eul_{x_0}(M;\mathbb R)$,
the co-Euler structure $[g,\alpha]\in\Eul^*_{\smash{x_0}}(M;\mathbb R)$ and the
cohomology class $[\omega]\in H^1(M;\mathbb R)$.
Indeed, in view of \fullref{L:S}\itemref{S:i}, 
\fullref{P:R}\itemref{P:R:o} and \itemref{P:R:i}
as well as $d\alpha=E(g)$, $\partial c=E(X)-\chi(M)x_0$ and
$\IND_\alpha(x_0)=\chi(M)x_0$ the quantity \eqref{rhoinv}
remains unchanged when $\omega$ is replaced by $\omega+dh$.
\fullref{P:R}\itemref{P:R:vii} and \eqref{E:rxxx}
show that it does not depend on the representative of the Euler 
structure $[X,c]$. In view of \fullref{P:R}\itemref{P:R:vi} 
and \eqref{E:rgg} it does not depend on the representative of the
co-Euler structure either. Thus \eqref{rhoinv} defines a coupling
$$
\mathbb T\co\Eul_{x_0}(M;\R)\times\Eul^*_{\smash{x_0}}(M;\R)\to H_1(M;\R).
$$
From the very definition we have
\begin{equation}\label{Tdef}
\bigl\langle[\omega],\mathbb T([X,c],[g,\alpha])\bigr\rangle
=\int_M\omega\wedge(X^*\Psi(g)-\alpha)-\int_c\omega,
\end{equation}
where $\omega$ is any representative of $[\omega]$ which vanishes locally
around the zeros of $X$ and vanishes locally around the base point $x_0$.
Moreover, we have
\begin{equation}\label{Tequi}
\mathbb T(\e+\sigma,\e^*+\beta)
=\mathbb T(\e,\e^*)-\sigma+\PD^{-1}(\beta)
\end{equation}
for all $\e\in\Eul_{x_0}(M;\R)$, $\e^*\in\Eul^*_{\smash{x_0}}(M;\R)$, 
$\sigma\in H_1(M;\R)$ and $\beta\in H^{n-1}(M;\mathcal O_M)$. 
Here $\PD\co H_1(M;\R) \to H^{n-1}(M;\mathcal O_M)$ denotes the 
Poincar\'e duality isomorphism.

We have the following affine version of Poincar\'e duality:

\begin{proposition}\label{affPD}
There is a natural isomorphism
$$
P\co\Eul_{x_0}(M;\R)\to\Eul^*_{\smash{x_0}}(M;\R)
$$
affine over the Poincar\'e duality isomorphism $\PD\co H_1(M;\R) \to H^{n-1}(M;\mathcal O_M)$.
More precisely, for $\sigma\in H_1(M;\R)$ and 
$\e\in\Eul_{x_0}(M;\R)$ we have 
$$
P(\e+\sigma)=P(\e)+\PD(\sigma).
$$
This affine isomorphism intertwines the involution on $\Eul_{x_0}(M;\R)$
with the involution on $\Eul^*_{\smash{x_0}}(M;\R)$, that is, $P(\nu(\e))=\nu(P(\e))$ for all
$\e\in\Eul_{x_0}(M;\R)$. Moreover, $\mathbb T(\e,\e^*)=P^{-1}(\e^*)-\e$.
\end{proposition}

\begin{proof}
In view of \eqref{Tequi} we may define $P$ by $\mathbb T(\e,P(\e))=0$.
The equivariance property and the last equation follow at once. The 
bijectivity of $P$ follows from the bijectivity of $\PD$.
The statement $P(\nu(\e))=\nu(P(\e))$ is a consequence of \eqref{E:nupsi};
see \fullref{R:caneul} and \fullref{R:cancoeul}.
\end{proof}

Let us collect the observations from the last two sections.

\begin{theorem}\label{T:eul}
Let $(M,x_0)$ be a closed connected base-pointed manifold.
\begin{enumerate}
\item\label{T:eul:iii}
Let $\pi_0(\vf(M,x_0))$ denote the set of connected components of the space
of vector fields which just vanish in $x_0$, equipped with the 
$C^\infty$--topology. If $\dim M>2$ then we have a bijection
$$
\pi_0(\vf(M,x_0))\to\Eul_{x_0}(M;\Z),\qquad [X]\mapsto [X,0].
$$
\item\label{T:eul:iv}
Let $\pi_0(\vf_0(M))$ denote the set of connected components of the space
of no\-where vanishing vector fields, equipped with the $C^\infty$--topology.
If $\chi(M)=0$ and $\dim M>2$ we have a surjection
$$
\pi_0(\vf_0(M))\to\Eul_{x_0}(M;\Z),\qquad [X]\mapsto[X,0].
$$
\item\label{T:eul:i}
There exists a natural isomorphism
$$
P\co\Eul_{x_0}(M;\R)\to\Eul^*_{\smash{x_0}}(M;\R)
$$
affine over the Poincar\'e duality 
$\PD\co H_1(M;\R) \to H^{n-1}(M;\mathcal O_M)$, that is,
$P(\e+\sigma)=P(\e)+\PD(\sigma)$, 
for all $\sigma\in H_1(M;\R)$.
\item\label{T:eul:ii}
The assignment 
$
\T\co\Eul_{x_0}(M;\R)\times\Eul^*_{\smash{x_0}}(M;\R)\to H_1(M;\R)
$ which is given by $\T(\e,\e^*)=P^{-1}(\e^*)-\e$
satisfies 
$$
\langle[\omega],\T(\e,\e^*)\rangle
=\int_M\omega\wedge(X^*\Psi(g)-\alpha)-\int_c\omega
$$
where $\e=[X,c]$, $\e^*=[g,\alpha]$
and $\omega\in\Omega^1(M;\mathbb R)$ is any representative of the class
$[\omega]\in H^1(M;\mathbb R)$ which vanishes locally around 
$x_0$ and the zeros of $X$.
\end{enumerate}
\end{theorem}

We continue to denote by $P$ resp.\ $\PD$ the composition of $P$ resp.\ 
$\PD$ with the canonical mapping $\Eul_{x_0}(M;\Z)\to\Eul_{x_0}(M;\R)$ resp.\ 
the homomorphism $H_1(M;\Z)\to H_1(M;\R)$. The image 
$P(\Eul_{x_0}(M;\Z))\subseteq\Eul_{x_0}(M;\R)$ is referred to as the 
\emph{lattice of integral co-Euler structures\/}. If two Euler structures 
$\e_1,\e_2\in\Eul_{x_0}(M;\Z)$ satisfy $P(\e_1)=P(\e_2)$ then 
$\e_2-\e_1\in H_1(M;\Z)$ is a torsion element.

\subsection{Milnor--Turaev torsion}\label{SS:MT-tor}

Let $M$ be a closed connected manifold of dimension $n$. Fix a base point
$x_0\in M$ and let 
$\Gamma=\pi_1(M,x_0)$ denote the fundamental group. Let $\mathbb K$
be a field of characteristic zero and let $V$ be a finite-dimensional vector space over
$\mathbb K$. Let $\rho\co\Gamma\to\GL(V)$ be a representation.
Consider the graded $\mathbb K$--vector space $H^*(M;\rho)$,
its determinant line $\det H^*(M;\rho)$, and the one-dimensional 
$\mathbb K$--vector space\footnote{For the definition of (graded) 
determinant lines see the footnote in the introduction.}
\begin{equation}\label{E:Detr}
\Det_{x_0}(M;\rho):=\det H^*(M;\rho)\otimes(\det V)^{-\chi(M)}. 
\end{equation}

\begin{remark}\label{R:acDet}
If $\chi(M)=0$ then 
$\Det_{x_0}(M;\rho)=\det H^*(M;\rho)$. If $H^*(M;\rho)=0$ then
$\Det_{x_0}(M;\rho)$ is canonically isomorphic to $\mathbb K$.
\end{remark}

Suppose we have a smooth triangulation $\tau$ of $M$. It gives rise to a 
cochain complex of $\mathbb K$--vector spaces $C^*_\tau(M;\rho)$ which 
computes the cohomology $H^*(M;\rho)$. By standard linear algebra we get a 
canonical isomorphism of $\mathbb K$--vector spaces
\begin{equation}\label{caniso}
\varphi^\rho_\tau\co\det C^*_\tau(M;\rho)
\to\det H(C^*_\tau(M;\rho))=\det H^*(M;\rho).
\end{equation}
This isomorphism is described in detail by Farber and Turaev
\cite[section 2, formula (2.2)]{FT99} where a sign
correction\footnote{first introduced in \cite{Tu86}} $(-1)^{N(C)}$ of 
previously used definitions is added; see also \eqref{E:rst} in the 
appendix of this paper.

Recall the Euler vector field $X_\tau$ of $\tau$ from
\fullref{SS:triang}. Its zero set $\mathcal X_\tau$ coincides
with the set of barycenters of $\tau$. For a cell 
$\sigma$ of $\tau$ we let $x_\sigma\in\mathcal X_\tau$ denote its barycenter.
We have $\IND_{X_\tau}(x_\sigma)=(-1)^{\dim\sigma}$ for every cell $\sigma$
of $\tau$. Suppose we have given an Euler structure $\e\in\Eul_{x_0}(M;\Z)$. 
By \fullref{vector-field}, see also Turaev \cite {Tu90},
one can choose a collection of paths $\pi_x$ from $x_0$ to $x\in\mathcal X_\tau$ so 
that with $c:=\sum_{x\in\mathcal X_\tau}\IND_{X_\tau}(x)\pi_x$ we have 
$\e=[X_\tau,c]$.

Let $v$ be a nonzero element in $\det V$.\footnote{Note 
that a frame (basis) in $V$ determines such an
element. Note that if instead of $v$ we pick up a
frame $\epsilon$ of $V$, we will obtain a base in
$C^*_\tau(M;\rho)$.} Using parallel transport along $\pi_x$ we get a
nonzero element in $\det C^*_\tau(M;\rho)$ with a sign ambiguity.
If the set of barycenters $\mathcal X_\tau$ were ordered we would get 
a well-defined nonzero element in $\det C^*_\tau(M;\rho)$.

Let $\mathfrak o$ be a cohomology orientation of $M$, \ie an
orientation of $\det H^*(M;\R)$. We say an ordering of the zeros is compatible
with $\mathfrak o$ if the nonzero element in $\det C^*_\tau(M;\R)$ provided
by this ordered base is compatible with the orientation $\mathfrak o$ via the
canonical isomorphism
$$
\varphi^\mathbb R_\tau\co\det C^*_\tau(M;\R)\to\det H(C^*_\tau(M;\R))=\det H^*(M;\R).
$$

So given a collection of paths $\pi_x$ with $\e=[X_\tau,c]$, 
$c=\sum_{x\in X_\tau}\IND_{X_\tau}(x)\pi_x$, an ordering of
$\mathcal X_\tau$ compatible with $\mathfrak o$, and a nonzero element $v\in\det V$ 
one obtains (using Milnor's construction modified by Turaev) a nonzero element in 
$\det C^*_\tau(M;\rho)$ corresponding to a nonzero element in 
$\det H^*(M;\rho)$ via \eqref{caniso}. We thus get a mapping
\begin{equation}\label{hommap}
\det V\setminus 0\to\det H^*(M;\rho)\setminus 0.
\end{equation}
This mapping is obviously homogeneous of degree $\chi(M)$. A straight forward 
calculation shows that it is independent  on the choice of $\pi_x$ 
or the ordering on $\mathcal X_\tau$ provided they define the same  
$\e$ and $\mathfrak o$. 
As a matter of fact this mapping does not depend on $\tau$ either but 
only on the Euler structure $\e$ and the cohomology orientation $\mathfrak o$.
In fact it depends on $\mathfrak o$ only when $\dim V$ 
is odd; cf \eqref{E:mto} below.
This is a nontrivial fact, and its proof is contained in Milnor \cite{M66} and 
Turaev \cite{Tu86} for the acyclic case, and implicit in the existing literature 
\cite{FT99}, but see also \fullref{C:app} in the appendix.

\begin{definition}[Milnor--Turaev torsion]\label{D:MT-tor}
The element in $\Det_{x_0}(\!M;\rho)$ corre\-spond\-ing to the homogeneous mapping \eqref{hommap} 
is referred to as the \emph{Milnor--Turaev torsion\/} and will be denoted by
$\tau_\comb^{\rho,\e,\mathfrak o}$; see \eqref{E:Detr}.
\end{definition}

We obviously have 
$$
\tau_\comb^{\rho,\e+\sigma,\mathfrak o}
=\tau_\comb^{\rho,\e,\mathfrak o}\cdot[\det\circ\rho](\sigma)^{-1},
\quad \text{for all } \sigma\in H_1(M;\Z).
$$
See \fullref{R:KT} for the sign. Here $[\det\circ\rho]\co H_1(M;\Z)\to\mathbb K^*$
denotes the homomorphism induced from the homomorphism
$\det\circ\rho\co\Gamma\to\mathbb K^*$. 
In particular, if $\rho$ is an unimodular 
representation then $\tau_\comb^{\rho,\e,\mathfrak o}$ is
independent of $\e$. Moreover
\begin{equation}\label{E:mto}
\tau_\comb^{\rho,\e,-\mathfrak o}=(-1)^{\dim V}\tau_\comb^{\rho,\e,\mathfrak o}.
\end{equation}

\begin{remark}\label{R:acy-MT}
Note that when $H^*(M;\rho)=0$ then, in view of the canonical isomorphism
$\Det_{x_0}(M;\rho)=\mathbb K$ of \fullref{R:acDet}, the Milnor--Turaev 
torsion is an element in $\mathbb K^*$. To stay with traditional notation we  
will denote the inverse of this number by $\mathcal T_\comb^{\e,\mathfrak o}(\rho)
:=(\tau_\comb^{\rho,\e,\mathfrak o})^{-1}$, and later consider it as a
function on the space of acyclic representations. 
\end{remark}

\subsection{Modified Ray--Singer metric}\label{SS:mRS-met}

Let us briefly recall the definition of the Ray--Singer torsion \cite{RS71}.
Suppose $\mathbb F=(F,\nabla)$ is a real or complex
vector bundle $F$ equipped with a flat connection $\nabla$. Let $g$ be a 
Riemannian metric on $M$ and let $\mu$ be a Hermitian metric on $F$.
The Ray--Singer torsion $T_\an(\nabla,g,\mu)$ is the positive
real number given by
$$
\log T_\an(\nabla,g,\mu):=\frac12\sum_k(-1)^{k+1}k\log{\det}'\Delta^k,
$$
where $\Delta^k$ is the Laplacian acting in degree $k$ of the elliptic complex 
$(\Omega^*(M;F),d_\nabla)$ equipped with the scalar product induced from 
the Riemannian metric $g$ and the Hermitian metric $\mu$, and 
${\det}'\Delta^k$ denotes its zeta regularized determinant, ignoring the
zero eigenvalues.

Recall that $H^*(M;\mathbb F)$ canonically identifies with the harmonic
forms and thus inherits a Hermitian metric from $g$ and $\mu$. Let us write
$\protect{\|\cdot\|_\Hodge^{\mathbb F,g,\mu}}$ for the induced Hermitian metric on
$\det H^*(M;\mathbb F)$.
The Ray--Singer metric is the Hermitian metric 
$$
\|\cdot\|_\RS^{\mathbb F,g,\mu}
:=T_\an(\nabla,g,\mu)\cdot\|\cdot\|_\Hodge^{\mathbb F,g,\mu}
$$
on $\det H^*(M;\mathbb F)$, see \cite{BZ92}.

Suppose $\mathbb K$ is $\R$ or $\C$. Let $V$ be a finite-dimensional vector
space over $\mathbb K$. Recall that every representation of the fundamental
group $\rho\co\Gamma\to\GL(V)$ induces a canonical vector bundle $F_\rho$
equipped with a canonical flat connection $\nabla_\rho$. They are obtained from the
trivial bundle $\tilde M\times V\to\tilde M$ and the trivial connection by
passing to the $\Gamma$ quotient spaces. Here $\tilde M$ is the canonical
universal covering provided by the base point $x_0$. The $\Gamma$--action is
the diagonal action of deck transformations on $\tilde M$ and the action 
$\rho$ on $V$. The fiber of $F_\rho$ over $x_0$ identifies canonically with
$V$. The holonomy representation determines a right $\Gamma$--action on the
fiber of $F_\rho$ over $x_0$, \ie an anti homomorphism $\Gamma\to\GL(V)$.
When composed with the inversion in $\GL(V)$ we get back the representation
$\rho$. The pair $(F_\rho, \nabla_\rho)$ will be denoted by $\mathbb F_\rho$.

A Hermitian metric on $F_\rho$ induces a Hermitian metric 
$\smash{\protect{\|\cdot\|_{\smash{x_0}}^\mu}}$ on $(\det V)^{-\chi(M)}$. Moreover,
on $\det H^*(M;\mathbb F_\rho)=\det H^*(M;\rho)$
we have the Ray--Singer metric $\protect{\|\cdot\|_\RS^{\smash{\mathbb F_{\rho},g,\mu}}}$.
Suppose $\e^*\in\Eul_{x_0}(M;\R)$ is a co-Euler structure. Choose
$\alpha$ so that $\e^*=[g,\alpha]$; see
\fullref{metric}. Recall the Kamber--Tondeur form 
$\omega(\nabla_\rho,\mu)\in\Omega^1(M;\R)$ from \fullref{SS:KT}.
Consider the Hermitian metric
\begin{equation}\label{E:mRS-met}
\|\cdot\|_\an^{\rho,\e^*}:=
\|\cdot\|_\RS^{\mathbb F_\rho,g,\mu}\otimes
\|\cdot\|_{x_0}^\mu
\cdot e^{-S(\omega(\nabla_\rho,\mu),\alpha)}
\end{equation}
on $\Det_{x_0}(M;\rho)$; see \eqref{E:Detr} and \fullref{L:S}.

\begin{lemma}\label{L:indep}
The Hermitian metric \eqref{E:mRS-met} only depends on $\rho$
and the co-Euler structure $\e^*=[g,\alpha]$ represented by $g$ and $\alpha$.
\end{lemma}

\begin{proof}
This follows from well-known anomaly formulas for the Ray--Singer metric 
\cite{BZ92}. We first keep $\mu$ fixed and consider two Riemannian metrics
$g_1$ and $g_2$ as well as $\alpha_1$ and $\alpha_2$ so that 
$[g_1,\alpha_1]=[g_2,\alpha_2]=\e^*$. In view of \fullref{L:S} this implies
\begin{equation}\label{E:712}
S(\omega,\alpha_2)-S(\omega,\alpha_1)=\int_M\omega\wedge\cs(g_1,g_2)
\end{equation}
for every closed one-form $\omega\in\Omega^*(M;\mathbb R)$.
In this situation \cite[Theorem~0.1]{BZ92} tells, with our notation,
$$
\log\left(\frac{\|\cdot\|_\RS^{\mathbb F_\rho,g_2,\mu}}
{\|\cdot\|_\RS^{\mathbb F_\rho,g_1,\mu}}\right)^2=
2\int_M\omega(\nabla_\rho,\mu)\wedge\cs(g_1,g_2).
$$
Together with \eqref{E:712} we obtain
\begin{multline*}
\log\frac{\|\cdot\|_\RS^{\mathbb F_\rho,g_2,\mu}\otimes
\|\cdot\|_{x_0}^\mu
\cdot e^{-S(\omega(\nabla_\rho,\mu),\alpha_2)}}
{\|\cdot\|_\RS^{\mathbb F_\rho,g_1,\mu}\otimes
\|\cdot\|_{x_0}^\mu
\cdot e^{-S(\omega(\nabla_\rho,\mu),\alpha_1)}}
\\=\int_M\omega(\nabla_\rho,\mu)\wedge\cs(g_1,g_2)
-S(\omega(\nabla_\rho,\mu),\alpha_2)+S(\omega(\nabla_\rho,\mu),\alpha_1)
=0.
\end{multline*}
Thus \eqref{E:mRS-met} does not depend on the representative of $\e^*$.

Next we keep $g$ and $\alpha$ fixed and consider two Hermitian metrics
$\mu_1$ and $\mu_2$. Then from \eqref{E:omm} we have
$\omega(\nabla_\rho,\mu_2)-\omega(\nabla_\rho,\mu_1)
=-\frac12d\log V(\mu_1,\mu_2)$ and \fullref{L:S} implies
\begin{multline}\label{E:abc}
S(\omega(\nabla_\rho,\mu_2),\alpha)-
S(\omega(\nabla_\rho,\mu_1),\alpha)
=-\frac12S(d\log V(\mu_1,\mu_2),\alpha)
\\
=-\frac12\chi(M)\log V(\mu_1,\mu_2)(x_0)+\frac12\int_M\log V(\mu_1,\mu_2)E(g).
\end{multline}
In this situation \cite[Theorem~0.1]{BZ92} states, with our notation,
$$
\log\left(\frac{\|\cdot\|_\RS^{\mathbb F_\rho,g,\mu_2}}
{\|\cdot\|_\RS^{\mathbb F_\rho,g,\mu_1}}\right)^2
=\int_M\log V(\mu_1,\mu_2)E(g).
$$
Together with \eqref{E:abc} and since 
$\bigl(\|\cdot\|^{\mu_2}_{x_0}\big/\|\cdot\|_{x_0}^{\mu_1}\bigr)^2
=\bigl(V(\mu_1,\mu_2)(x_0)\bigr)^{-\chi(M)}$
we obtain
\begin{multline*}
\log\left(\frac{\|\cdot\|_\RS^{\mathbb F_\rho,g,\mu_2}\otimes
\|\cdot\|_{x_0}^{\mu_2}
\cdot e^{-S(\omega(\nabla_\rho,\mu_2),\alpha)}}
{\|\cdot\|_\RS^{\mathbb F_\rho,g,\mu_1}\otimes
\|\cdot\|_{x_0}^{\mu_1}
\cdot e^{-S(\omega(\nabla_\rho,\mu_1),\alpha)}}\right)^2
=\int_M\log V(\mu_1,\mu_2)E(g)-
\\-\chi(M)\log V(\mu_1,\mu_2)(x_0)
-2S(\omega(\nabla_\rho,\mu_2),\alpha)
+2S(\omega(\nabla_\rho,\mu_1),\alpha)
=0.
\end{multline*}
Hence \eqref{E:mRS-met} is independent of $\mu$.
\end{proof}

For $\beta\in H^{n-1}(M;\Or_M)$ we have
$$
\|\cdot\|_\an^{\rho,\e^*+\beta}
=\|\cdot\|_\an^{\rho,\e^*}
\cdot e^{\langle\log|[\det\circ\rho]|,\PD^{-1}(\beta)\rangle}
$$
(\fullref{L:S}\itemref{S:i}) where $\log|[\det\circ\rho]|\in H^1(M;\R)$ is the Kamber--Tondeur class
corresponding to the homomorphism $\log|[\det\circ\rho]|\co H_1(M;\R)\to\R$ as in
\fullref{R:KT}.

\begin{definition}[Modified Ray--Singer metric]\label{D:mRS-met}
The Hermitian metric $\protect{\|\cdot\|_\an^{\rho,\e^*}}$ on $\Det_{x_0}(M;\rho)$ will be referred 
to as the \emph{modified Ray--Singer metric\/}; see \eqref{E:Detr} and \eqref{E:mRS-met}.
\end{definition}

\begin{remark}
Suppose $n$ is odd and recall the canonical co-Euler structure
$\e^*_{\smash{\can}}$ in $\Eul^*(M;\mathbb R)$. Clearly, in this situation
$\|\cdot\|_\an^{\rho,\e^{\smash{*}}_\can}=\|\cdot\|_\RS^{\smash{\mathbb F_\rho},g,\mu}$. From \fullref{L:indep} we recover the well-known 
fact, that $\|\cdot\|^{\smash{\mathbb F_\rho},g,\mu}_\RS$ is 
independent of $g$ and $\mu$, provided $n$ is odd.
\end{remark}

The Bismut--Zhang theorem \cite[Theorem~0.2]{BZ92} can now be reformulated as
a statement comparing two topological invariants, one defined with the help
of analysis, and one defined with the help of a triangulation.

\begin{theorem}[Bismut--Zhang]\label{T:BZ}
Let $(M,x_0)$ be a closed connected manifold with base point, and 
let $\rho\co\pi_1(M,x_0)\to\GL(V)$ be a finite-dimensional real or complex
representation. Suppose $\e\in\Eul_{x_0}(M;\Z)$ is an Euler
structure and $\e^*\in\Eul^*_{\smash{x_0}}(M;\R)$ a co-Euler structure.  
Then one has
$$
\|\tau_\comb^{\rho,\e,\mathfrak o}\|_\an^{\rho,\e^*}
=
e^{\langle\log|[\det\circ\rho]|,\T(\e,\e^*)\rangle}.
$$
Particularly, if $P(\e)=\e^*$ then
$\protect{\|\tau_\comb^{\rho,\e,\mathfrak o}}\|_\an^{\rho,\e^*}=1$.
\end{theorem}

For an alternative proof of the (original) Bismut--Zhang theorem; see 
also Burghelea, Friedlander and Kappeler \cite{BFK01}.

\begin{remark}\label{R:acRS}
If $H^*(M;\rho)=0$ then
$\Det_{x_0}(M;\rho)=\mathbb K$; see \fullref{R:acDet}. In this case 
the modified Ray--Singer metric gives rise to a real number 
$\mathcal T_\an^{\e^*}(\rho):=\|1\|_\an^{\rho,\e^*}>0$.
From the very definition we have
$$
\mathcal T_\an^{\e^*}(\rho)=
T_\an(\nabla_\rho,g,\mu)\cdot e^{-S(\omega(\nabla_\rho,\mu),\alpha)}
$$
where $\nabla_\rho$ denotes the flat connection on the associated bundle
$F_\rho$, $g$ is a Riemannian metric on $M$, $\mu$ is a
Hermitian structure on $F_\rho$, and $\alpha$ is such that $[g,\alpha]=\e^*$.
Recall the number $\mathcal T_\comb^{\e,\mathfrak o}(\rho)\in\mathbb K^*$ from
\fullref{R:acy-MT}. Assuming $P(\e)=\e^*$ the Bismut--Zhang theorem tells
$$
1=\|\tau^{\rho,\e,\mathfrak o}_\comb\|^{\rho,\e^*}_\an
=|\mathcal T_\comb^{\e,\mathfrak o}(\rho)|^{-1}\|1\|^{\rho,\e^*}_\an
=|\mathcal T_\comb^{\e,\mathfrak o}(\rho)|^{-1}\mathcal T_\an^{\e^*}(\rho),
$$ 
that is $|\mathcal T_\comb^{\e,\mathfrak o}(\rho)|=\mathcal T_\an^{\e^*}(\rho)$.
\end{remark}

\section{Torsion and representations}\label{S:rep}

{\bf Warning}\qua 
Following standard conventions in algebraic geometry \cite{Ha}
we reserve the name \emph{affine algebraic variety\/} for an irreducible 
variety in the affine space, and the term \emph{algebraic set\/} for any finite 
union of varieties. By a rational function on an algebraic set we simply mean 
a rational function on each irreducible component.

\subsection{Complex representations}\label{SS:rep}

Suppose $\Gamma$ is a finitely presented group with generators $g_1,\dotsc,g_r$ 
and relations $R_i(g_1,g_2,\dotsc,g_r)=e$, $i=1,\dotsc,p$, and $V$ be a
complex vector space of dimension $N$. Let $\Rep(\Gamma;V)$ be the set of 
linear representations of $\Gamma$ on $V$, \ie group homomorphisms 
$\rho\co\Gamma\to\GL_\C(V)$. By identifying $V$ to $\C^N$ this set is in a
natural way an algebraic set inside the space $\C^{\smash{rN^2+1}}$ given
by $pN^2+1$ equations. Precisely if $A_1,\dotsc,A_r,z$ represent the 
coordinates in $\C^{\smash{rN^2+1}}$ with $A:=(a^{ij})$, $a^{ij}\in\C$, so 
$A\in\C^{\smash{N^2}}$ and $z\in\C$, then the equations defining $\Rep(\Gamma;V)$ are 
\begin{eqnarray*}
z\cdot\det(A_1)\cdot\det(A_2)\cdots\det(A_r)&=&1
\\
R_i(A_1,\dotsc,A_r)&=&\id,\qquad i=1,\dotsc,p
\end{eqnarray*}
with each of the equalities $R_i$ representing $N^2$ polynomial equations.

Suppose $\Gamma=\pi_1(M,x_0)$, $M$ a closed manifold.
Denote by $\Rep^M_0(\Gamma;V)$  the set of representations $\rho$ 
with $H^*(M;\rho)=0$ and notice that they form a Zariski
open set in $\Rep(\Gamma;V)$. Denote the Zariski closure of this set by  
$\Rep^M(\Gamma;V)$. This is again an algebraic set 
which depends only on the homotopy type of $M$, 
can be empty (for example if $\chi(M)\neq0$) but in some interesting
cases can be the full space $\Rep(\Gamma;V)$; cf \fullref{SS:maptor}.

\begin{definition}[$V$--acyclicity]\label{D:V-ac}
Given a complex vector space $V$  we say that the closed
manifold $M$ is \emph{$V$--acyclic\/} if there exist representations 
$\rho\in\Rep(\Gamma;V)$ with $H^*(M;\rho)=0$.
Equivalently $\Rep^M(\Gamma;V)\neq\emptyset$.
\end{definition}

\subsection{The space of cochain complexes and Turaev functions}\label{SS:Tfunc}

Let $(k_0, k_1,\dotsc,k_n)$ be a string of nonnegative integers which satisfy
the following requirements:
\begin{eqnarray}
k_0-k_1+k_2\mp\cdots+(-1)^nk_n&=&0
\label{RE:0}
\\
k_i-k_{i-1}+ k_{i-2}\mp\cdots+(-1)^ik_0&\geq&0
\quad\text{for any $i\leq n-1$.}
\label{RE:1}
\end{eqnarray}

Denote by $\mathbb D(k_0,\dotsc,k_n)$
the collection of cochain complexes of the form
$$
C=(C^*,d^*):
0\to C^0 \xrightarrow{d^0} 
C^1\xrightarrow{d^1}
\cdots 
\xrightarrow{d^{n-2}}
C^{n-1}
\xrightarrow{d^{n-1}} C^n\to 0
$$
with $C^i:=\C^{k_i}$,
and by $\mathbb D_\ac(k_0,\dotsc,k_n)\subseteq\mathbb D(k_0,\dotsc,k_n)$ 
the subset of acyclic complexes. The cochain complex $C$ is determined by 
the collection $\{d^i\}$ of linear maps $d^i\co\C^{k_i}\to\C^{k_{i+1}}$.
If $\mathbb D(k_0,\dotsc,k_n)$ is regarded as the subset of 
$\{d^i\}\in\mathbb L(k_0,\dotsc,k_n)
=\bigoplus_{i=0}^{n-1}L(\C^{k_i}, \C^{k_{i+1}})$\footnote{We denote 
by $L(V,W)$ the space of linear maps from $V$ to $W$.}
which satisfy the quadratic equations 
$$
d^{i+1}\cdot d^i=0,
$$
then $\mathbb D(k_0,\dotsc,k_n)$ is an affine algebraic 
set given by degree two homogeneous polynomials and 
$\mathbb D_\ac(k_0,\dotsc,k_n)$ is a Zariski open set
in $\mathbb D(k_0,\dotsc,k_n)$.
We denote its Zariski closure by $\hat{\mathbb D}_\ac(k_0,\dotsc,k_n)$.

\begin{proposition}\label{P:irred}
\
\begin{enumerate}
\item
$\mathbb D_\ac(k_0,\dotsc,k_n)$ is a connected smooth quasialgebraic set of dimension 
$$
k_0\cdot k_1 + (k_1- k_0)\cdot k_2+\cdots +(k_{n-1}-k_{n-2} +\cdots \pm k_0)\cdot k_{n}.
$$
\item
$\hat {\mathbb D}_\ac (k_0,\dotsc, k_n)$ is an irreducible algebraic
set, hence an affine algebraic variety 
\end{enumerate}
\end{proposition}

\begin{proof}
The map 
$$ 
\pi_0\co\mathbb D_\ac (k_0,\dotsc, k_n)\to\Emb(C^0,C^1)
$$
which associates to $C\in\mathbb D_\ac(k_0,\dotsc,k_n)$
the linear map $d^0$, is a bundle whose fiber is isomorphic to 
$\mathbb D_\ac (k_1-k_0,k_2,\dotsc,k_n)$.
Note that $\Emb(C^0,C^1)$ is connected, smooth and of dimension 
$k_0\cdot k_1$. By induction this implies that $\mathbb
D_\ac(k_0,\dotsc,k_n)$ is a connected smooth quasialgebraic set and the
dimension is as claimed. This shows the first part. The second part follows
from the well-known fact \cite[page~21]{GH} that a complex algebraic set, 
whose nonsingular part is connected with respect to standard topology, 
must be irreducible.
\end{proof}

For any cochain complex in $C\in \mathbb D_\ac (k_0,\dotsc,k_n)$ denote by 
$$
B^i:=\img(d^{i-1})\subseteq C^i= \C^{k_i}
$$
and consider the short exact sequence 
$$
0\to B^i \xrightarrow{\textrm{inc}} C^i \xrightarrow {d^i} B^{i+1}\to 0.
$$ 
Choose a bases $b_i$ for each $B_i$, and choose lifts $\overline{b}_{i+1}$ of
$b_{i+1}$ in $C^i$ using $d^i$, \ie $d^i(\overline{b}_{i+1})= b_{i+1}$. 
Clearly $\{b_i,\overline b_{i+1}\}$ is a base of $C^i$.

Consider the base $\{b_i,\overline{b}_{i+1}\}$ as a collection of vectors in 
$C^i= \C^{k_i}$ and write them as a columns of a matrix $[b_i,\overline{b}_{i+1}]$. 
Define the torsion of the acyclic complex $C$, by 
$$
\mathfrak t(C):=(-1)^N\prod^n_{i=0}\det[b_i,\overline{b}_{i+1}]^{(-1)^i}
$$
where $(-1)^N$ is Turaev's sign; see Farber and Turaev \cite{FT99} or appendix~\ref{S:app}.
The result is independent of the choice of 
the bases $b_i$ and of the lifts $\overline{b}_i$ \cite{M66,FT99} 
and leads to the function 
\begin{equation}\label{E:007}
\mathfrak t\co\mathbb D_\ac(k_0,\dotsc,k_n)\to\C\setminus 0.
\end{equation} 
Turaev provided a simple formula for this function \cite{Tu01} which  
permits to recognize $\mathfrak t$ as the restriction of a rational function on  
$\hat{\mathbb D}_\ac(k_0,\dotsc,k_n)$.

Following Turaev \cite{Tu01}, one considers symbols
$\alpha\equiv(\alpha_0,\alpha_1,\dotsc,\alpha_n)$ 
with $\alpha_i$ subsets the set of integers $\{1,2,\dotsc,k_i\}$ 
so that
\begin{eqnarray*}
\alpha_n&=&\emptyset
\\ 
\sharp(\alpha_i)&=&k_{i+1}-\sharp(\alpha_{i+1}).
\end{eqnarray*}
In Turaev \cite{Tu90} such an $\alpha$ is called a $\tau$--chain. 
There are, of course, only finitely many 
$\tau$--chains for a given collection $(k_0,\dotsc,k_n)$ as above.

For any cochain complex $C\in{\mathbb D}(k_0,\dotsc, k_n)$, and in fact for any element 
$d\equiv\{d^i\}$ in $\mathbb L(k_0,\dotsc,k_n)$, a $\tau$--chain defines the collection
of square matrices $A_i(d^i)$ representing the $\sharp(\alpha_i)\times
(k_{i+1}-\sharp(\alpha_{i+1}))$--minor of the matrix $d^i$ 
consisting of the columns whose indices are in $\alpha_i$ and
the rows not in $\alpha_{i+1}$.

For $d\equiv\{d^i\}\in\mathbb L(k_0,\dotsc,k_n)$
one considers the expression 
$$
F_{\alpha}(d)=\epsilon(\alpha)\prod_{i=0}^n(\det A_i(d^i))^{(-1)^{i+1}}
$$
where $\epsilon(\alpha)=\pm 1$ and is chosen according to the strategy explained 
by Turaev \cite[Remark~2.4, page 9]{Tu01}.  
Clearly $F_\alpha$ is a rational function on $\mathbb L(k_0,\dotsc,k_n)$.

The cochain complex is called $\alpha$--nondegenerate if all quantities 
$\det A_i(d^i)$ are non\-zero. Turaev has shown that:
\begin{enumerate}
\item\label{T:i}
For an $\alpha$--nondegenerate cochain complex 
$C\in \mathbb D_\ac(k_0,\dotsc,k_n)$, 
$\mathfrak t(C)= F_\alpha(d)$ where $d=\{d_i\}$ is collection of the differentials in the 
complex $C$ \cite[Theorem~2.2]{Tu01}.
\item\label{T:ii}
If $U_\alpha$ denotes the subsets of $\mathbb D_\ac (k_0,\dotsc, k_n)$ consisting
$\alpha$--nondegenerate complexes then $U_\alpha$ is open, and 
$\bigcup_{\alpha} U_\alpha= \mathbb D_\ac (k_0,\dotsc, k_n)$. 
\item\label{T:iii}
For any $\tau$--chain $\alpha$ there exists a cochain complex 
$C\in\mathbb D_\ac(k_0,\dotsc,k_n)$ which is
$\alpha$--nondegenerate.
\end{enumerate}

Using \itemref{T:i}--\itemref{T:iii} and the irreducibility of 
$\hat{\mathbb D}_\ac(k_0,\dotsc,k_n)$ 
we conclude: every $F_\alpha$ restricts to a rational function on 
$\hat{\mathbb D}_\ac(k_0,\dotsc,k_n)$; as rational
functions, these restrictions agree, are regular on $\mathbb D_\ac(k_0,\dotsc,k_n)$,
and they coincide with \eqref{E:007} on $\mathbb D_\ac(k_0,\dotsc,k_n)$.
Therefore \eqref{E:007} is the restriction of a rational function
on $\hat{\mathbb D}_\ac(k_0,\dotsc,k_n)$ with zeros and poles contained
in $\hat{\mathbb D}_\ac(k_0,\dotsc,k_n)\setminus\mathbb D_\ac(k_0,\dotsc,k_n)$.

Consider now a smooth triangulation $\tau$ of $M$ whose set of simplices of dimension 
$q$ is denoted by $\mathcal X_q$ and the collection of integers 
$k_i=\sharp(\mathcal X_i)\cdot\dim V$ where $V$ is a fixed complex vector space.
Let $\e\in\Eul_{x_0}(M;\Z)$ be an Euler structure, 
$\mathfrak o$ a cohomology orientation, and suppose that $M$ is $V$--acyclic.
Then the integers $(k_0,\dotsc,k_n)$ satisfy \eqref{RE:0} and \eqref{RE:1}. 
As in \fullref{SS:MT-tor}, choose a collection of paths
$\pi_{\e}\equiv\{\pi_x\}$ from $x_0$ to the barycenters $x$, 
so that $\e=[X_\tau,c]$ where $c$ is the Euler chain defined by
$\{\pi_x\}$. Choose also an ordering $o$ of the barycenters (zeros of $X_\tau$) 
compatible with $\mathfrak o$, and a framing $\epsilon$ of $V$.

Consider the chain complex $(C^*_\tau(M;\rho),d_\tau(\rho))$ associated with
the triangulation $\tau$ which computes the cohomology $H^*(M;\rho)$.
Using $\pi_x$, $o$ and $\epsilon$ one can identify $C^q_\tau(M;\rho)$ with
$\C^{k_q}$. We obtain in this way a map 
$$ 
t_{\pi_{\e},o,\epsilon}\co\Rep(\Gamma;V)\to\mathbb D(k_0,\dotsc,k_n)
$$
which sends $\Rep^M(\Gamma;V)\setminus\Sigma(M)$ to $\mathbb D_\ac (k_0,\dotsc,k_n)$. 
Here $\Sigma(M)\subseteq\Rep(\Gamma;V)$ denotes the Zariski closed subset of 
representations $\rho$ which are not local minima for the function
$\rho\mapsto\sum_{i=0}^n\dim H^i(M;\rho)$.
A look at the explicit definition of $d_\tau(\rho)$ implies 
that $t_{\pi_{\e},o,\epsilon}$ is actually a regular map between two algebraic sets.
Change of $\epsilon$, $\pi_{\e}$ and $o$ changes the map $t_{\pi_{\e},o,\epsilon}$. 
Since we are in the case $H^*(M;\rho)=0$ we easily recognize that the
composition $\mathfrak t\cdot t_{\pi_{\e},o,\epsilon}$, when restricted to 
$\Rep^M(\Gamma;V)\setminus \Sigma(M)$, is exactly the Milnor--Turaev 
torsion $\mathcal T_\comb^{\e,\mathfrak o}$ (\fullref{R:acy-MT})
hence independent of all these choices $\pi_{\e},o,\epsilon$.

We summarize these observations in the following theorem.

\begin{theorem}\label{T:ratio}
Let $M$ be a closed manifold, $x_0$ a base point, $\Gamma:=\pi_1(M,x_0)$ the
fundamental group, $\e\in\Eul_{x_0}(M;\Z)$ an Euler structure,
$\mathfrak o$ a cohomology orientation, and let $V$ be a finite-dimensional 
complex vector space. Consider the Milnor--Turaev torsion
$\mathcal T_\comb^{\e,\mathfrak o}\co\Rep^M(\Gamma;V)\setminus\Sigma(M)\to\C\setminus 0$.
\medskip\begin{enumerate}
\item\label{T:ratio:i}
$\mathcal T_\comb^{\e,\mathfrak o}$ is a rational
function on the algebraic set $\Rep^M(\Gamma;V)$ whose poles and 
zeros are contained in $\Sigma(M)$.
\item\label{T:ratio:ii}
The absolute value of
$\mathcal T_\comb^{\e,\mathfrak o}$ is the modified Ray--Singer torsion 
$\mathcal T_\an^{\e^*}$, where $\e^*=P(\e)\in\Eul_{\smash{x_0}}^*(M;\R)$; see \fullref{R:acRS}. More precisely, for
a representation $\rho\in\Rep^M(\Gamma;V)\setminus\Sigma(M)$ one has
$$
|\mathcal T_\comb^{\e,\mathfrak o}(\rho)|
=\mathcal T_\an^{\e^*}(\rho)
=T_\an(\nabla_\rho,g,\mu,\alpha)
=T_\an(\nabla_\rho,g,\mu)\cdot e^{-S(\omega(\nabla_\rho,\mu),\alpha)}
$$
where $(F_\rho,\nabla_\rho)$ denotes the flat complex vector bundle associated to the
representation $\rho$, $g$ is a Riemannian metric on
$M$, $\mu$ is a Hermitian structure on $F_\rho$, $\omega(\nabla_\rho,\mu)$ is the
Kamber--Tondeur form, and $\alpha$ is such that the pair $(g,\alpha)$
represents the co-Euler structure $\e^*=P(\e)$.
\item\label{T:ratio:iii}
For $\sigma\in H_1(M;\Z)$, and $\rho\in\Rep^M(\Gamma;V)\setminus\Sigma(M)$ we have
$$\phantom{999}
\mathcal T_\comb^{\e+\sigma,\mathfrak o}(\rho)
=\mathcal T_\comb^{\e,\mathfrak o}(\rho)
\cdot[\det\circ\rho](\sigma)
\quad\text{and}\quad
\mathcal T_\comb^{\e,-\mathfrak o}(\rho)=
(-1)^{\dim V}\cdot\mathcal T_\comb^{\e,\mathfrak o}(\rho)
$$
where $[\det\circ\rho]\co H_1(M;\Z)\to\C\setminus 0$ is the homomorphism
induced from the map $\det\circ\rho\co\Gamma\to\C\setminus 0$.
\end{enumerate}
\end{theorem}

\subsection{Milnor--Turaev torsion for mapping tori}\label{SS:maptor}

Let $N$ be a closed connected manifold and $\varphi\co N\to
N$ a diffeomorphism. Define the mapping torus $M=N_\varphi$ by gluing the
boundaries of $N\times I$ with the help of $\varphi$, more precisely 
identifying $(x,1)$ with $(\varphi(x),0)$. The manifold $M$ comes with a
projection $p\co M\to S^1$, and an embedding $N\to M$, $x\mapsto(x,0)$. 
Note that $\pi_1(S^1,0)= \Z$.
Let $\K$ be a field of characteristics zero, and suppose
$V$ is a finite-dimensional vector space over $\K$.
For $\rho\in\Rep(\Z;V)=\GL(V)$ define
$$
P^k_\varphi(\rho):=\det\Bigl(
\varphi^*\circ\rho_*-\id\co H^k(N;V)\to H^k(N;V)
\Bigr)
$$
and the Lefschetz zeta function
$$
\zeta_\varphi(\rho)
:=\frac{\prod_{\text{$k$ even}}P^k_{\varphi}(\rho)}
{\prod_{\text{$k$ odd}}P^k_{\varphi}(\rho)}.
$$
This is a rational function $\zeta_\varphi\co\GL(V)\to\K$.

Choose a base point $x_0\in M$ with $p(x_0)=0\in S^1$. Every representation
$\rho\in\Rep(\Z; V)$ gives rise to a representation $p^*\rho$ of 
$\pi_1(M,x_0)$ on $V$. It is not hard to see that
$H^*(M;p^*\rho)=0$ if and only if $P^k_\varphi(\rho)\neq0$ for all $k$.
Indeed, we have a long exact sequence, the Wang sequence:
\begin{equation}\label{E:Wang}
\cdots\to H^*(\!M;\!p^*\rho)\to H^*(\!N;\!V)
\xrightarrow{\varphi^*\circ\rho_*-\id} 
H^*(\!N;\!V)
\to H^{*+1}(\!M;\!p^*\rho)\to\cdots 
\end{equation}
Recall that this sequence is derived from the long exact sequence
of the pair $(M,N)$, and the isomorphism
$H^*(M,N;p^*\rho)=H^*(N\times I,N\times\partial I;V)
=H^{*-1}(N;V)$. The latter is the Thom isomorphism with the standard
orientation on $I$, and the first comes from a trivialization of
the coefficient system on $N\times I$ from $N\times\{0\}$ upwards.

Let $\omega:=p^*dt\in\Omega^1(M)$. The space of vector fields 
$X$ on $M$ satisfying $\omega(X)<0$ is contractible and hence 
defines an Euler structure $\e\in\Eul(M;\Z)$; cf 
\fullref{T:eul}\itemref{T:eul:iv}.

Recall that a finite-dimensional long exact sequence of $\K$ vector spaces 
$$
\cdots\to A^n\to B^n\to C^n\to A^{n+1}\to\cdots
$$ 
induces a canonical isomorphism of graded determinant lines
\begin{equation}\label{E:123}
u\co\det(A^*)\otimes\det(C^*)\to\det(B^*),
\qquad
a\otimes b\to u(a\otimes b).
\end{equation}
Indeed, if one regards this sequence as an acyclic  
cochain complex $$0\to D^0\to D^1\cdots\to D^n\to D^{n+1}\to\cdots$$
one obtains a canonical isomorphism $\det (D^*)\to\mathbb K$; see \eqref{E:rst} in the appendix.
(Since the complex is acyclic the Turaev sign correction does not appear.)
In view of the identification of $\det(D^\ast)$ with 
$\det (A^\ast)\otimes (\det (B^\ast))^{-1} \otimes \det(C^\ast)$  
one obtains the canonical isomorphism \eqref{E:123}.  
If in addition $B^*=C^*$ we obtain a canonical isomorphism 
$\det (A^*)\to\K$, $a\to\tr(u(a\otimes\cdot))$. 
Applying this to the Wang exact sequence for the trivial one-dimensional real representation
$$ 
\cdots\to H^*(M;\R)\to H^*(N;\R)\to H^*(N;\R)\to H^{*+1}(M;\R)\to \cdots
$$ 
we obtain a canonical isomorphism $\det H^*(M;\R)=\R$ hence a canonical cohomology 
orientation $\mathfrak o$.

The rest of this section is dedicated to the proof of the following proposition;
see also Fried \cite{Fr87}.

\begin{proposition}\label{P:maptor}
Let $\e$ and $\mathfrak o$ denote the canonical Euler structure and
cohomology orientation on the mapping torus $M$ introduced above.
If $H^*(M;p^*\rho)=0$, then 
$\mathcal T_\comb^{\e,\mathfrak o}(p^*\rho)
=(-1)^{\dim(V)\cdot z_\varphi}\zeta_\varphi(\rho)$ where
$$
z_\varphi=\sum_q\dim H^q(N;\R)\cdot
\dim\ker\bigl(\varphi^*-\id\co H^q(N;\R)\to H^q(N;\R)\bigr).
$$
\end{proposition}

To get a hold on the Euler structure $\e$ we start with the following lemma. 
\begin{lemma}\label{L:eul}
Suppose $\tilde\tau$ is a smooth triangulation of $M$ such that $N$
is a subcomplex. For every simplex $\sigma\in\tilde\tau$ let $\pi_\sigma$
be a path in $N\times[0,1)$ considered as path in $M$ 
from the base point $x_0\in N$ to the 
barycenter $x_\sigma$. Consider the Euler structure
$\tilde\e=[X_{\tilde\tau},c_{\tilde\tau}]$ where 
$c_{\tilde\tau}:=\sum_{\sigma\in\tilde\tau}(-1)^{\dim\sigma}\pi_\sigma$.
If $\e$ is the Euler structure described above then 
$[\det\circ p^*\rho](\tilde\e-\e)=1$. Particularly,
$\mathcal T_\comb^{\e,\mathfrak o}(p^*\rho)
=\mathcal T_\comb^{\tilde\e,\mathfrak o}(p^*\rho)$. The same statement remains true when 
$\tilde\tau$ is a smooth cell structure as in \fullref{R:vect}.
\end{lemma}

\begin{proof}
Let $X$ be a vector field on $M$ with $\omega(X)<0$ and recall that
$\e=[X,0]$.
Since $N$ is a subcomplex we find a sufficiently small $\epsilon>0$ such that
$N\times[1-\epsilon,1)\subseteq M$ does not contain barycenters of simplices of 
$\tilde\tau$. In view of P2 in \fullref{SS:triang} we may also assume
$\omega(X_{\tilde\tau})<0$ on $N\times[1-\epsilon,1)$. Let us write
$N_\epsilon:=N\times\{1-\epsilon\}\subseteq M$.
We thus find a homotopy $\X$ from 
$X$ to $X_{\tilde\tau}$ which has the property that $\omega(\X_t)<0$ on 
$N_\epsilon$ for all $t\in I$. We conclude that $\cs(X,X_{\tilde\tau})$ has a
representative supported
in $M\setminus N_\epsilon$. Since $\tilde\e-\e\in H_1(M;\Z)$ is
represented by $c_{\tilde\tau}-\cs(X,X_{\tilde\tau})$ it too has a
representative supported in $M\setminus N_\epsilon$. We conclude 
$p_*(\tilde\e-\e)=0\in H_1(S^1;\Z)$, and hence 
$[\det\circ p^*\rho](\tilde\e-\e)=1$. The last assertion follows from
\fullref{T:ratio}\itemref{T:ratio:iii}.
\end{proof}

\begin{lemma}\label{L:dan}
Let $\varphi_i\co N\to N$ be diffeomorphisms, let $\e_i$ resp.\ $\mathfrak o_i$
denote the Euler structure resp.\ cohomology orientation
on $N_{\varphi_i}$ introduced above, and suppose
$H^*(N_{\varphi_i};p^*\rho)=0$, $i=a,b$. Then
$$
\mathcal T_\comb^{N_{\varphi_b},\e_b,\mathfrak o_b}(p^*\rho)\cdot
\mathcal T_\comb^{N_{\varphi_a},\e_a,\mathfrak o_a}(p^*\rho)^{-1}
=(-1)^{\dim(V)\cdot(z_{\varphi_b}-z_{\varphi_a})}
\zeta_{\varphi_b}(\rho)\cdot\zeta_{\varphi_a}(\rho)^{-1}.
$$
\end{lemma}

\begin{proof}
We may assume, by changing $\varphi_i$ up to isotopy, that there are smooth
triangulations $\tau_0$ and $\tau_1$ of $N$ such that both $\varphi_a$ and
$\varphi_b$ are simplicial. Let $\tau$ be a smooth triangulation
on $N\times[0,1]$ such that $N\times\{0\}$ and $N\times\{1\}$ are
subcomplexes and the triangulation induced by $\tau$ on $N\times\{0\}$
resp.\ $N\times\{1\}$ is exactly $\tau_0$ resp.\ $\tau_1$.
Let $\tilde\tau_a$ and $\tilde\tau_b$ denote the \emph{cell structures\/} on
$N_{\varphi_a}$ and $N_{\varphi_b}$ induced by $\tau$.\footnote{Observe 
that if $K$ is a simplicial complex with $L_1$ and $L_2$ two disjoint 
subcomplexes and $\varphi\co L_1\to L_2$ a simplicial map the quotient space 
$K_\varphi$ obtained by identifying $x\in L_1$ to $\varphi(x)\in L_2$ 
is a cell complex whose cells are in bijective correspondence to the 
open simplices of $K\setminus L_1$.}  
Note that $\tilde\tau_a$ and $\tilde\tau_b$ are smooth cell structures 
(of the type considered in \fullref{R:vect}) but 
not smooth triangulations. The characteristic maps are smooth 
embeddings of compact affine simplices; they satisfy the axioms of 
a cell complex but not of a triangulation. Also recall that the torsion 
$\mathcal T_\comb^{N_{\varphi_i},\e,\mathfrak o}(p^*\rho)$ can be computed using
any cell structure \cite{M66}. A simple inspection of the Turaev sign 
corrections indicates that such correction can be performed  
mutatis-mutandis in this more general case.

We denote by 
$\mathcal X$, $\mathcal X_0$ and $\mathcal X_1$ the set of simplices (cells)
of $\tau$, $\tau_0$ and $\tau_1$. 
The cells of both $\tilde\tau_a$ and $\tilde\tau_b$ are indexed by 
$\mathcal X\setminus\mathcal X_0$.
Hence the graded vector spaces
$C^*_{\tilde\tau_a}(N_{\varphi_a};p^*\rho)$ and
$C^*_{\tilde\tau_b}(N_{\varphi_b};p^*\rho)$ are tautologically the same.  
Moreover the relative cochain complexes 
$C^*_{\tilde\tau_a}(N_{\varphi_a},N;p^*\rho)$,
$C^*_{\tilde\tau_b}(N_{\varphi_b},N;p^*\rho)$ and
$C^*_\tau(N\times I,N\times\partial I;p^*\rho)$ are tautologically the same;  
their cells are indexed by $\mathcal X\setminus\{\mathcal X_0\cup\mathcal
X_1\}$ and the differentials coincide.

The two short exact sequences of cochain complexes
$$
C^*_\tau(N\times I,N\times\partial I;p^*\rho)
= C^*_{\tilde\tau_i}(N_{\varphi_i},N;p^*\rho)
\to C^*_{\tilde\tau_i}(N_{\varphi_i};p^*\rho)\to
C^*_{\tau_1}(N;V)
$$
$i=a,b$, give rise to commuting fusion diagrams as explained in 
\fullref{L:detses} in the appendix. These make up the upper half
of the following diagram:
$$\disablesubscriptcorrection
\scalebox{0.77}{
\xymatrix{
\det C^*_{\tilde\tau_a}(N_{\varphi_a};p^*\rho)
\ar[d]_-{\varphi_{\tilde\tau_a}}
&
\det C^*_\tau(N\times I,N\times\partial I;p^*\rho)\otimes\det C^*_{\tau_1}(N;V) 
\ar[d]^-{\varphi\otimes\varphi_{\tau_1}}
\ar[l]_-{(-1)^{y_a}\psi_a}
\ar[r]^-{(-1)^{y_b}\psi_b}
&
\det C^*_{\tilde\tau_b}(N_{\varphi_b};p^*\rho)
\ar[d]^-{\varphi_{\tilde{\tau}_b}}
\\
\K=\det H^*(N_{\varphi_a};p^*\rho)
&
\det H^*(N\times I,N\times\partial I;p^*\rho)\otimes\det H^*(N;V) 
\ar[l]_-{\varphi_{\mathcal H_a^*}}
\ar[r]^-{\varphi_{\mathcal H_b^*}}
&
\det H^*(N_{\varphi_b};p^*\rho)=\K
\\
& 
\det H^{*-1}(N;V)\otimes\det H^*(N;V)=\mathbb K
\ar[ul]^-{\zeta_{\varphi_a}(\rho)^{-1}}
\ar[ur]_-{\zeta_{\varphi_b}(\rho)^{-1}}
\ar@{=}[u]
}}
$$
Note that \fullref{L:detses} is here applied to short exact sequences
where the middle complex $C_1^*$ is acyclic. Therefore
\fullref{L:detses} also implies
\begin{equation}\label{E345}
{y_a}\equiv{y_b}\mod2.
\end{equation}
The vertical identification in the lower part of the diagram
is explained in the text below \eqref{E:Wang}. Because of \eqref{E:Wang}
the lower triangles commute too.

For every simplex $\sigma\in\mathcal X\setminus\mathcal X_0$ let 
$\pi_\sigma$ be a path in $N\times[0,1)$ from the base point $x_0\in N$ 
to the barycenter $x_\sigma$ and consider the Euler chains 
$c_{\tilde\tau_a}$ resp.\ $c_{\tilde\tau_b}$
defined by the collections of paths $\{\pi_\sigma\}$ viewed in 
$N_{\varphi_a}$ resp.\ $N_{\varphi_b}$ (by the formula
$c_{\tilde\tau}:=\sum_{\sigma\in\tilde\tau}(-1)^{\dim\sigma}\pi_\sigma$.)
Consider the diagram 
\begin{equation}\disablesubscriptcorrection\label{D:aa}
\scalebox{0.74}{\xymatrix{
{\det C^*_{\tilde\tau_a}(N_{\varphi_a};p^*\rho)\phantom{1}}
&
{\phantom{1}\det C^*_\tau(N\times I,N\times\partial I;p^*\rho)\otimes\det C^*_{\tau_1}(N;V)\phantom{1}} 
\ar[l]_-{(-1)^{y_a}\psi_a}
\ar[r]^-{(-1)^{y_b}\psi_b}
&
{\phantom{1}\det C^*_{\tilde\tau_b}(N_{\varphi_b};p^*\rho)}
\\
\K=\det V^{\chi}
\ar[u]_-{\beta_{\varphi_a}}&
\K=\det V^{\chi}
\ar[l]_-{(-1)^{\dim(V)\cdot z_a}}
\ar[r]^-{(-1)^{\dim(V)\cdot z_b}}
&
\det V^{\chi}=\K
\ar[u]_-{\beta_{\varphi_b}}
}}
\end{equation}
where $\beta_{\varphi_a}$ and $\beta_{\varphi_b}$ are induced from the  
homogeneous maps of degree $\chi=\chi(N_{\varphi_a})=\chi(N_{\varphi_b})=0$ 
which associates to a frame of $V$ the
equivalence class of bases induced from the frame, the collection of paths 
$\{\pi_\sigma\}$ and the ordering of the simplices in 
$\mathcal X\setminus \mathcal X_0$ representing the canonical cohomology 
orientation $\mathfrak o_a$ resp.\ $\mathfrak o_b$ of 
$N_{\varphi_a}$ resp.\ $N_{\varphi_b}$.

We would like to prove that the diagram is commutative.
If so, in view of  \fullref{L:eul}, \fullref{D:MT-tor}  
and \fullref{R:acy-MT}, by putting together the two diagrams we 
conclude the result as stated.

It is easy to  note the commutativity  up to sign (independent on 
the representation $\rho$) of the second diagram, since the Euler 
structures $\e_a$ resp.\ $\e_b$ are realized by the same Euler 
chain given by $\{\pi_\sigma\}$; cf \fullref{L:eul}. 
We then conclude the result up to a sign independent of the 
representation $\rho$.

To decide about the sign we have to compare the orderings of $\mathcal
X\setminus\mathcal X_0$ induced from the canonical cohomology orientation
on $N_{\varphi_a}$ and on $N_{\varphi_b}$. To do this we consider
the fusion diagram for the trivial one-dimensional real representation:
$$
\disablesubscriptcorrection
\scalebox{0.9}{\xymatrix{
\det C^*_{\tilde\tau_a}(N_{\varphi_a};\R)
\ar[d]_-{\varphi_{\tilde\tau_a}^\R}
&
\det C^*_\tau(N\times I,N\times\partial I;\R)\otimes\det C^*_{\tau_1}(N;\R) 
\ar[d]^-{\varphi^\R\otimes\varphi_{\tau_1}^\R}
\ar[l]_-{(-1)^{y_a^\R}\psi^\R_a}
\ar[r]^-{(-1)^{y_b^\R}\psi^\R_b}
&
\det C^*_{\tilde\tau_b}(N_{\varphi_b};\R)
\ar[d]^-{\varphi_{\tilde{\tau}_b}^\R}
\\
\det H^*(N_{\varphi_a};\R)
&
\det H^*(N\times I,N\times\partial I;\R)\otimes\det H^*(N;\R) 
\ar[l]_-{\varphi^\R_{\mathcal H_a^*}}
\ar[r]^-{\varphi^\R_{\mathcal H_b^*}}
&
\det H^*(N_{\varphi_b};\R)
\\
&\det H^{*-1}(N;\mathbb R)\otimes\det H^*(N;\mathbb R)=\mathbb R
\ar@{=}[u]
}}
$$
Let $o_i$ denote the equivalence class of orderings of $\mathcal
X\setminus\mathcal X_0$
induced by the cohomology orientation $\mathfrak o_i$,
$i=a,b$. The commutativity of the previous diagram tells
\begin{equation}\label{E:567}
\sign(o_b/o_a)=(-)^{y_b^{\mathbb R}-y_a^{\mathbb R}}.
\end{equation}
Diagram \eqref{D:aa} will commute if and only if
$$
(-1)^{\dim(V)\cdot(z_b-z_a)}(-1)^{y_b-y_a}\sign(o_b/o_a)^{\dim(V)}=1.
$$
In view of \eqref{E345} and \eqref{E:567} it thus suffices to show
\begin{equation}\label{E:O}
y_b^\R-y_a^\R\equiv z_b-z_a\mod2.
\end{equation}
To establish this relation let us write for $i=a,b$
\begin{align}
\notag
x_i^q&:=\dim\ker\bigl((\varphi_i^*-\id)\co H^q(N;\mathbb R)\to 
H^q(N;\mathbb R)\bigr)
\\ \label{E:A}
f^q_{i,1}&:=\dim\img\bigl(H^q(C^*_\tau(N\times I,N\times\partial I;\mathbb R))
\to H^q(C^*_{\tilde\tau_i}(N_{\varphi_i};\mathbb R))\bigr)
\\ \label{E:B}
f_{i,2}^q&:=\dim\img\bigl(H^q(C^*_{\tilde\tau_i}(N_{\varphi_i};\mathbb R))\to
H^q(C^*_{\tau_1}(N;\mathbb R))\bigr)
\end{align}
where the homomorphisms in \eqref{E:A} and \eqref{E:B} are the ones
from the long exact cohomology sequence associated with the short exact
sequence of complexes
$$
0\to C^*_\tau(N\times I,N\times\partial I;\mathbb R)
\to C^*_{\tilde\tau_i}(N_{\varphi_i};\mathbb R)\to
C^*_{\tau_1}(N;\mathbb R)\to0.
$$
Note that for $(\varphi_i^*-\id)\co H^q(N;\mathbb R)\to H^q(N;\mathbb R)$,
we have
$$
\dim\ker (\varphi_i^*-\id)
\\=\dim\coker (\varphi_i^*-\id),
$$
and thus \eqref{E:Wang} gives
\begin{equation}\label{E:C}
f^q_{i,2}=x^q_i\quad\text{and}\quad f^q_{i,1}=x^{q-1}_i,\qquad i=a,b.
\end{equation}
Let us also introduce
\begin{align*}
b_0^q&:=\dim B^q\bigl(C^*_\tau(N\times I,N\times\partial I;\mathbb R)\bigr)
\\
b_2^q&:=\dim B^q\bigl(C^*_{\tau_1}(N;\mathbb R)\bigr)
\end{align*}
where $B^q$ denotes the space of boundaries in degree $q$.
From \fullref{L:detses} and \eqref{E:C} we obtain
\begin{align}
\notag
y^{\mathbb R}_b-y^{\mathbb R}_a
&=N\bigl(C^*_{\tilde\tau_b}(N_{\varphi_b};\mathbb R)\bigr)
-N\bigl(C^*_{\tilde\tau_a}(N_{\varphi_a};\mathbb R)\bigr)
\\\notag&\qquad
+\sum_q(f^q_{b,1}-f^q_{a,1})b_2^q+\sum_q(f^q_{b,2}-f^q_{b,1})b_0^{q+1}
\\\label{E:D}
&=N\bigl(C^*_{\tilde\tau_b}(N_{\varphi_b};\mathbb R)\bigr)
-N\bigl(C^*_{\tilde\tau_a}(N_{\varphi_a};\mathbb R)\bigr)
+\sum_q(x^{q-1}_b-x^{q-1}_a)(b_2^q+b_0^q).
\end{align}
Let us further introduce
\begin{align*}
\alpha_0^q&:=\sum_{j\geq q}\dim C^j_\tau(N\times I,N\times\partial I;\mathbb R)
\\
\alpha_{i,1}^q&:=\sum_{j\geq q}\dim C^j_{\tilde\tau_i}(N_{\varphi_i};\mathbb R)
\qquad i=a,b
\\
\alpha_2^q&:=\sum_{j\geq q}\dim C^j_{\tau_1}(N;\mathbb R)
\end{align*}
and note that
$
    \alpha^q_{a,1}=\alpha^q_{b,1}=\alpha^q_0+\alpha^q_2.
$
Finally set
\begin{align*}
\beta_0^q&:=\sum_{j\geq q}\dim H^j\bigl(C^*_\tau(N\times I,N\times\partial I;\mathbb R)\bigr)
\\
\beta^q_{i,1}&:=\sum_{j\geq q}\dim
H^j\bigl(C^*_{\tilde\tau_i}(N_{\varphi_i};\mathbb R)\bigr)
\qquad i=a,b
\\
\beta^q_2&:=\sum_{j\geq q}\dim H^j\bigl(C^*_{\tau_1}(N;\mathbb R)\bigr).
\end{align*}
Note that $\beta^q_0=\beta_2^{q-1}$ and
$\beta_2^{q-1}-\beta^q_2=\dim H^{q-1}(N;\mathbb R)$
and hence
\begin{equation}\label{E:F}
\beta^q_0+\beta^q_2\equiv
\dim H^{q-1}(N;\mathbb R)\mod2.
\end{equation}
From \eqref{E:Wang} we find $\dim
H^j(C^*_{\tilde\tau_i}(N_{\varphi_i};\mathbb R))=x^j_i+x^{j-1}_i$, $i=a,b$,
and thus
\begin{equation*}
\beta^q_{i,1}\equiv x^{q-1}_i\mod 2\qquad i=a,b.
\end{equation*}
From this and $\alpha^q_{a,1}=\alpha^q_{b,1}=\alpha^q_0+\alpha^q_2$ we obtain
\begin{align*}
N\bigl(C^*_{\tilde\tau_b}(N_{\varphi_b};\mathbb R)\bigr)
-N\bigl(C^*_{\tilde\tau_a}(N_{\varphi_a};\mathbb R)\bigr)
&=\sum_{q\geq0}\alpha_{b,1}^q\beta^q_{b,1}
-\sum_{q\geq0}\alpha^q_{a,1}\beta^q_{a,1}
\\
&\equiv\sum_q(\alpha_0^q+\alpha_2^q)(x^{q-1}_b-x^{q-1}_a)\mod 2.
\end{align*}
Combining this with \eqref{E:D} we obtain
\begin{equation}\label{E:H}
y^{\mathbb R}_b-y^{\mathbb R}_a
\equiv\sum_q(\alpha^q_0+b^q_0+\alpha^q_2+\beta^q_2)(x^{q-1}_b-x^{q-1}_a)
\mod2.
\end{equation}
In view of the fact that for a finite cochain complex of 
finite-dimensional vector spaces $(C^*,d^*)$ the total rank of $C^*$ is
the same mod 2 as the total rank of its cohomology we have
$$
\alpha^q_0+b^q_0\equiv\beta^q_0\mod2
\quad\text{and}\quad
\alpha^q_2+b^q_2\equiv\beta^q_2\mod2
$$
and with the help of \eqref{E:F} we find
$$
\alpha^q_0+b^q_0+\alpha^q_2+b^q_2\equiv\beta^q_0+\beta^q_2
\equiv\dim H^{q-1}(N;\mathbb R)\mod2.
$$
Together with \eqref{E:H} this implies \eqref{E:O}, and the proof is completed.
\end{proof}

\begin{lemma}
If $H^*(N_\id;p^*\rho)=0$, then
$\mathcal T_\comb^{N_\id,\e,\mathfrak o}(p^*\rho)
=(-1)^{\dim(V)\cdot z_\id}\zeta_\id(\rho)$.
\end{lemma}

\begin{proof}
Equip the interval $I$ with the three cell structure. Let $\tau_1$ be a 
smooth triangulation of $N$, $\tau$ the product cell 
structure on $N\times I$ and $\tilde\tau$ the induced (product)
cell structure on $N_\id=N\times S^1$. Then the complex 
$C^*_\tau(N\times I,N\times\partial I;p^*\rho)$ is the suspension of 
$C^*_{\tau_1}(N;p^*\rho)$. Using \fullref{L:detses} and \fullref{L:eul}
it is easy to see that $\mathcal T_\comb^{\smash{N_\id},\e,\mathfrak o}(p^*\rho)=
\pm\zeta_\id(\rho)$.

To decide about the sign let $(-1)^{y^\rho}$
denote the sign in the fusion diagram of the short exact sequence
$$
0\to C^*_\tau(N\times I,N\times\partial I;p^*\rho)
\to C^*_{\tilde\tau}(N_\id;p^*\rho)
\to C^*_{\tau_1}(N;V)\to0
$$ 
and let $(-1)^{y^{\mathbb R}}$ denote the sign appearing in the
fusion diagram for the short exact sequence
$$
0\to C^*_\tau(N\times I,N\times\partial I;\mathbb R)
\to C^*_{\tilde\tau}(N_\id;\mathbb R)
\to C^*_{\tau_1}(N;\mathbb R)\to0.
$$
Similar to the proof of \fullref{L:dan} the statement about the 
sign follows once we have shown
\begin{equation}\label{E:AA}
\dim(V)\cdot y^\R-y^\rho\equiv\dim(V)\cdot z_\id\mod2.
\end{equation} 
Note that 
\begin{equation}\label{E:BB}
z_\id=\sum_q(\dim H^q(N;\mathbb R))^2
\equiv\sum_q\dim H^q(N;\mathbb R)
\equiv\chi(N)\mod2.
\end{equation}
Moreover let $(-1)^{y^V}$ denote the sign in the fusion diagram of
\begin{equation}\label{E:DD}
0\to C^*_\tau(N\times I,N\times\partial I;V)
\to C^*_{\tilde\tau}(N_\id;V)
\to C^*_{\tau_1}(N;V)\to0.
\end{equation}
Note that $y^V\equiv\dim(V)y^{\mathbb R}\mod2$. 
In view of \eqref{E:BB} it thus suffices to show
\begin{equation}\label{E:CC}
y^V-y^\rho\equiv\dim(V)\chi(N)\mod2.
\end{equation}
Using the acyclicity of $C^*_{\tilde\tau_\id}(N_\id;p^*\rho)$ and the simple
nature of the long exact sequence associated to \eqref{E:DD} we 
obtain from \fullref{L:detses} 
\begin{align}
\notag
y^V-y^\rho
&=N(C^*_{\tilde\tau_\id}(N_\id;V))
+\sum_q\dim H^{q-1}(N;V)b_2^q
+\sum_q\dim H^q(N;V)b^{q+1}_0
\\&\label{E:EE}
=N(C^*_{\tilde\tau_\id}(N_\id;V))
+\sum_q\dim H^{q-1}(N;V)(b_2^q+b_0^q)
\\
&b_0^q:=\dim B^q\bigl(C^*_\tau(N{\times} I,N{\times}\partial I;p^*\rho))
=\dim B^q(C^*_\tau(N{\times} I,N{\times}\partial I;V)\bigr)
\tag*{\hbox{where}}
\\
&b_2^q:=\dim B^q\bigl(C^*_{\tau_1}(N;V)\bigr).
\notag
\end{align}
It is not hard to verify
\begin{align}\label{E:FF}
&N(C^*_{\tilde\tau_\id}(N_\id;V))\equiv
\sum_q(\alpha_0^q+\alpha_2^q)\dim H^{q-1}(N;V)\mod2
\\
\tag*{\hbox{where}}
\alpha_0^q:=
&\sum_{j\geq q}\dim C^j_\tau(N\times I,N\times\partial I;p^*\rho)
=\sum_{j\geq q}\dim C^j_\tau(N\times I,N\times\partial I;V)
\\ \notag
\alpha^q_2:=&\sum_{j\geq q}\dim C^j_{\tau_1}(N;V).
\end{align}
Let us introduce
\begin{align*}
\beta^q_0&:=
\sum_{j\geq q}\dim H^j\bigl(C^*_\tau(N\times I,N\times\partial I;p^*\rho)\bigr)
\\&\qquad=\sum_{j\geq q}\dim H^j\bigl(C^*_\tau(N\times I,N\times\partial I;V)\bigr)
\\
\beta^q_2&:=\sum_{j\geq q}\dim H^j\bigl(C^*_{\tau_1}(N;V)\bigr)
\end{align*}
and recall that $\alpha^q_0+b_0^q\equiv\beta^q_0\mod2$ as well as
$\alpha^q_2+b^q_2\equiv\beta^q_0\mod2$. 
From \eqref{E:EE} and \eqref{E:FF} and
the fact that $\beta^q_2+\beta^q_0=\beta^q_2+\beta^{q-1}_2\equiv
\dim H^{q-1}(N;V)\mod2$ we obtain
\begin{align*}
y^V-y^\rho&\equiv\sum_q\dim H^{q-1}(N;V)(\beta^q_2+\beta^q_0)\mod2
\\&\equiv\sum_q(\dim H^{q-1}(N;V))^2\mod2
\\&\equiv\sum_q\dim H^{q-1}(N;V)\mod2
\\&\equiv\chi(N)\dim(V)\mod2.
\end{align*}
Hence we have established \eqref{E:CC} and the proof is complete.
\end{proof}

\fullref{P:maptor} follows from the previous two lemmas, and the
simple observation that $H^*(N_\varphi;p^*\rho)=0$ implies 
$H^*(N_\id;p^*\rho)=0$. Indeed, $H^0(N_\varphi;p^*\rho)=0$ and
\eqref{E:Wang} imply that $\rho(1)\in\GL(V)$ does not have $1$ as an 
eigenvalue and then \eqref{E:Wang} implies that $H^*(N_\id;p^*\rho)=0$.

\begin{remark}
Let $y\in N$ and suppose $\varphi (y)=y$. Denote by 
$\alpha\co\pi_1(N,y)\to\pi_1(N,y)$ the isomorphism  induced by $\varphi$. 
Let $\rho_0\co\pi_1(N,y)\to\GL(V)$ be a representation and $A\in\GL(V)$ so 
that $\rho_0(\alpha (g), A(v))= A(\rho_0 (g, v))$, $g\in\pi_1(N,y)$, 
$v\in V$. Denote by $\varphi^\ast_{\tilde\rho}\co 
H^\ast(N;\rho_0)\to H^\ast(N;\rho_0)$ resp.\ 
$\tilde\rho\co\pi_1(N,y)\times_\alpha\Z\to\GL(V)$ the  isomorphism 
induced by  the triple $(\varphi, \rho_0, A)$ resp.\ the representation 
induced by the triple $(\alpha, \rho_0, A)$.
The arguments in this section lead to a  result  similar to 
\fullref{P:maptor} where $H^\ast(N;V)$ is replaced by 
$H^\ast(N;\rho_0)$ and $H^\ast(M; p^\ast\rho)$ by $H^\ast(M;\tilde\rho)$.
\end{remark}

\subsection[An invariant with values in R/pi Z]{An invariant with values in $\R/\pi\Z$}\label{SS:inv}

If $M$ is a closed manifold with $\Gamma=\pi_1(M,x_0)$ and $\rho$ an
acyclic (\ie $H^*(M;\rho)=0$) unimodular representation 
then $\mathcal T_\comb^{\e,\mathfrak o}$ is by 
\fullref{T:ratio}\itemref{T:ratio:iii} independent of $\e$ and its argument mod 
$\pi$ independent of $\mathfrak o$. Denote by 
$$
S(M;\rho):=\arg(\mathcal T_\comb^{\e,\mathfrak o})\in\R/\pi\Z.
$$
This is a topological invariant computable with the help of a triangulation.
\fullref{P:maptor} gives a simple formula for
the case $M$ is a mapping torus and $\rho$ is the pullback by $p$ of a
unimodular representation of $\Z$, and implicitly establishes its
nontriviality.
In a forthcoming paper an analytic construction of this invariant
will be given with the help of spectral geometry of nonpositive Laplacians.

More generally, if $\e$ is an Euler structure and 
$\rho_1,\rho_2\in\Rep^M_0(\Gamma;V)$ then 
$$
\arg(\mathcal T_\comb^{\e,\mathfrak o})(\rho_2)
-\arg(\mathcal T_\comb^{\e,\mathfrak o})(\rho_1)
\in\R/2\pi\Z
$$
is independent of $\mathfrak o$ and defines an invariant 
$S^\e(M;\rho_1,\rho_2)\in\R/2\pi\Z$ which can be derived from the 
function $\mathcal T_\an^{\e^*}$, $P(\e)=\e^*$,
and calculated with the help of any holomorphic path $\tilde\rho(z)$ between 
$\rho_1$ and $\rho_2$ as explain in the introduction,
\begin{equation}\label{E:ph}
S^{\e}(M; \rho_1, \rho_2)
=\Re\Biggl(2/{\mathbf i}\int_1^2\frac{\partial(\mathcal T_\an^{\e^*}\circ\tilde\rho)}
{\mathcal T_\an^{\e^*}\circ\tilde\rho}\Biggr)
\end{equation}
where $\partial\varphi$ denotes the one-form $\smash{\frac{\partial\varphi}{\partial z}}dz$.
As already pointed out in introduction  given $\rho_1$  and $\rho_2$ two representations 
even in the same connected component of
$\Rep^M_0(\Gamma;V)$ a holomorphic path might not exist, however indirectly one 
can benefit from such formula.

It is interesting to compare this invariant with the one defined using the
spectral flow in \cite{APS75} for an odd-dimensional manifold and a unitary
representation \cite{F00}.   
More will be said in the forthcoming paper; see also the recent preprints
\cite{BK05a,BK05,BK05b,BK06}.

\subsection{Marcsik's theorem}\label{SS:marcsik}

In this section we will use the Bismut--Zhang theorem and
\fullref{P:maptor} to give a proof of a theorem due to J\,Marcsik.
We continue to use the notation from \fullref{SS:maptor}.

Let $g$ be a Riemannian metric on $M$ and $\omega$ a real-valued closed one-form 
on $M.$ For $z\in\C$ denote by  
$d_\omega(z)\alpha:=d\alpha+z\omega\wedge\alpha$ the Witten deformed 
differential on $\Omega^*(M;\C)$ and by $\delta_\omega(z)$ its adjoint. 
Denote by
$$
\Delta_\omega(z):=d_\omega(z)\delta_\omega(z)+\delta_\omega(z)d_\omega(z)
$$ 
the Laplacian, and by $\Delta^k_\omega(z)$ the Laplacian on $k$ forms. 
Introduce the Ray--Singer torsion
$$
T_\an(\omega,g)(z):=
\exp\Bigl(\frac12\sum_k(-1)^{k+1}k\log\det{}'\Delta^k_\omega(z)\Bigr).
$$
Here ${\det}'\Delta_\omega^k(z)$ denotes the regularized determinant of the nonnegative
self-adjoint elliptic operator $\Delta_\omega^k(z)$ obtained by ignoring the zero modes.

Identifying $\C^*=\GL(\C)$ we have $P^k_\varphi(w)$ and $\zeta_\varphi(w)$ defined 
for all $w\in\C^*$; see \fullref{SS:maptor}.

\begin{theorem}[Marcsik \cite{M98}]\label{T:marcsik}
Let $N$ be a closed manifold, $\varphi\co N\to N$ a diffeomorphism, 
$M$ the mapping torus of $\varphi$, $p\co M\to S^1$ the canonical 
projection, and $\omega$ the pullback of the canonical volume form on $S^1$.
Let $g$ be a Riemannian metric on $M$, and suppose
$X$ is a vector field on $M$ with $\omega(X)<0$.
Then
$$
\log T_\an(\omega,g)(z)
=\log|\zeta_\varphi(e^z)|
+\Re(z)\int_M\omega\wedge X^*\Psi(g)
$$
for all $z\in\C$ which satisfy $P^k_\varphi(e^z)\neq 0$ for all 
$k$.\footnote{This statement is slightly
more precise than the one formulated in \cite{M98}.}
\end{theorem}

\begin{proof}
Consider the trivial line bundle $M\times\C$ over $M$. Let $\nabla^0$
denote the trivial flat connection, and let $\mu$ denote the standard
Hermitian structure on $M\times\C$. For $z\in\C$ define a connection
$\nabla^z:=\nabla+z\omega$ on $M\times\C$. This connection is flat for
$\omega$ is closed. Clearly $d^{\nabla^z}=d_\omega(z)$, and hence
$$
T_\an(\omega,g)(z)
=T_\an(\nabla^z,g,\mu).
$$
Define a co-Euler structure $\e^*:=[g,X^*\Psi(g)]\in\Eul^*(M;\R)$.
Recall the Euler structure $\e=[X,0]\in\Eul(M;\Z)$. 
From \fullref{T:eul}\itemref{T:eul:ii} we see that $\T(\e,\e^*)=0$.
Assume $P^k_\varphi(e^z)\neq0$ for all $k$. Then
the Bismut--Zhang theorem implies
$$
|\mathcal T_\comb^{\e,\mathfrak o}(e^z)|=
T_\an(\nabla^z,g,\mu)\cdot
e^{-S(\omega(\nabla^z,\mu),X^*\Psi(g))}.
$$
See \fullref{R:acRS}.
Note that $\nabla^z\mu=-2\Re(z\omega)\mu$, hence
$\omega(\nabla^z,\mu)=\Re(z)\omega$, and thus
$$
S(\omega(\nabla^z,\mu),X^*\Psi(g))=\Re(z)\int_M\omega\wedge X^*\Psi(g).
$$
From \fullref{P:maptor} we get 
$\mathcal T_\comb^{\e,\mathfrak o}(e^z)=\pm\zeta_\varphi(e^z)$.
Putting everything together and taking the logarithm completes the proof.
\end{proof}

\subsection{An application}\label{SS:apl}

The concept of holomorphic (or $\mathbb C$--differentiable)   
maps between open sets of complex Banach/Fr\'echet 
spaces of possibly infinite dimension is the same as in   
finite dimension and so are the concepts of meromorphic map, complex 
analytic set, complex analytic manifold \cite{E66,Ra,D,HP,KM}. 
The reader unfamiliar with the infinite-dimensional holomorphy should be aware that:

(i)\qua If $U\subseteq\C$ is open and $E$ is a Fr\'echet space, then
all possible definitions for holomorphic curves $U\to E$ 
yield the same concept; see for instance Kriegl and Michor \cite[Theorem~7.4]{KM}.
A curve $c\co U\to E$ is holomorphic if and only if $\partial c/\partial\bar z=0$
and this is the case if and only if  $\varphi\circ c\co U\to\C$ is holomorphic for every
continuous linear functional $\varphi\in E^*$.

(ii)\qua If $V\subseteq F$ is an open subset of a Fr\'echet space, and $E$
is another Fr\'echet space, then a map $f\co V\to E$ is holomorphic if and only if
it maps holomorphic curves in $V$ to holomorphic curves in $E$. In view of
(i) a mapping $f\co V\to E$ is holomorphic if and only if  for every continuous linear functional
$\varphi\in E^*$ and every holomorphic curve $c\co U\to V$ the composition
$\varphi\circ f\circ c$ is holomorphic. This is the
case if and only if, $f$ is continuous and holomorphic along affine complex lines in
$F$; see \cite[Theorem~7.19]{KM}, \cite[Chapter III sections~2 and 3]{HP}
or \cite[pages 57, 58]{D}.

(iii)\qua A subset $X\subseteq E$ of a Fr\'echet space is called complex analytic subset
if for every point $x\in X$ there exists an open neighborhood $U_x$ of $x$
in $E$, a Fr\'echet space $F_x$ and a holomorphic map $f_x\co U_x\to F_x$ such that
$X\cap U_x=f_x^{-1}(0)$ \cite[page 383]{D}.

(iv)\qua Suppose $X$ is a complex analytic subset of a Fr\'echet space $E$
and $Y$ is a complex analytic subset of a Fr\'echet space $F$.
A map $f\co X\to Y$ is called holomorphic if for every point $x\in X$
there exists an open neighborhood $U_x$ of $x$ in $E$ and a holomorphic
function $f_x\co U_x\to F$ whose restriction to $X$ coincides with $f$.

(v)\qua Suppose $X$ is a complex analytic subset of a Fr\'echet space $E$. 
$X$ is called locally irreducible if for any $x\in X$ the local ring of 
germs of holomorphic functions $\mathfrak o_x(X)$ is an integral domain.
This  concept is used mostly for finite-dimensional complex analytic sets.

If $S\subseteq X$ is a complex analytic subset of $X$ it is called proper for
$X$, if the closure of $X\setminus S$ is dense in $X$.

(vi)\qua If $X$ is a finite-dimensional irreducible complex analytic set  of a
Fr\'echet space $E$ and
$S\subseteq X$ a proper  analytic subset, a meromorphic  function on $X$
with poles in $S$ is a function $f\co X\setminus S\to \mathbb C$
s.t.\ for any $x\in X$ there exists a neighborhood $U_x$ and two holomorphic
functions $h_x\co U_x\to \C$ and
$g_x\co U_x\to \C$, $g_x\neq0$ on $U_x\setminus S$, so that 
$f|_{U_x\setminus S}=(h_x/ g_x )|_{U_x\setminus S}$.

If $X$ is a complex analytic subset of a Fr\'echet space $E$ and
$S\subseteq X$ a proper analytic subset, a meromorphic function on $X$
with poles in $S$ is a function $f\co X\setminus S\to\mathbb C$
which, when restricted to any finite-dimensional locally irreducible
analytic subset $V\subseteq X$ with $V\cap S$ proper in $V$, is  
meromorphic with poles in $V\cap S$.

The above definitions can be formulated  in the case where $E$ is replaced by a 
complex analytic Frech\'et manifold.

Suppose that $p\co F\to M$ is a smooth complex vector bundle which admits flat connections.
Denote by $\mathcal C(F)$ the space of all
connections on $F$. This is an affine Fr\'echet space over $\C$ 
when equipped with the $C^\infty$--topology. In view of (ii) above 
the map which associates to each connection its curvature is a holomorphic map from the 
affine Fr\'echet space $\mathcal C(F)$ to the Fr\'echet space
$\Omega^2(M;\End(F))$. 
Hence the set of flat connections $\mathcal F \mathcal C(F)$
is a closed complex analytic subset of $\mathcal C(F)$.

Let $x_0\in M$ be a base point, let $\Gamma:=\pi_1(M,x_0)$ denote the
fundamental group and let $V:=F_{x_0}$ denote the fiber over $x_0$.
The holonomy map $\mathcal F\mathcal C(F)\to\Rep(\Gamma;V)$ which associates to a flat 
connection $\nabla$ its holonomy 
representation $\rho_\nabla$, is a holomorphic map between complex analytic
sets. To see this we first choose $r$ smooth paths  
$\gamma_i\co [0,1]\to M$, $i=1,\dotsc,r$, with $\gamma_i(0)=\gamma_i(1)=x_0,$ which 
represent a collection of generators of $\Gamma=\pi_1(M,x_0)$.
For a connection $\nabla$ denote the isomorphisms 
defined by parallel transport along $\gamma_i$
by $P_i^\nabla\in\End(V)$, $V=F_{x_0}$. Define 
$\pi\co\mathcal C(F)\to\End(V)^r\times\mathbb C$ by 
$$
\pi(\nabla):=\Bigl(P^\nabla_1,\dotsc,P^\nabla_r, 
\det(P^\nabla_1)^{-1}\cdots\det(P^\nabla_r)^{-1}\Bigr).
$$ 
This is a holomorphic map in view of (ii) above. Since its restriction to
$\mathcal F\mathcal C(F)$ coincides with the holonomy map $\mathcal
F\mathcal C(F)\to\Rep(\Gamma;V)$, the latter is holomorphic.

Denote by $\Sigma(M)$ the set of representations $\rho\in\Rep(\Gamma;V)$ which 
are not local minima for the function $d(\rho):=\sum_{i=0}^n\dim H^i(M;\rho)$. 
Observe that this is an algebraic subset of $\Rep(\Gamma;V)$. In particular
it is a closed analytic subset of $\End(V)^r\times\mathbb C$. 
Let $\Sigma (F)$ denote the subset of flat connections $\nabla$
which are not a local minimum for the function 
$d\co\mathcal F\mathcal C(F)\to\Z$ defined by 
$d(\nabla):=\sum _{i=0}^n \dim H^i(M;\rho_{\nabla})$. 
Since $\Sigma(F)$ is the intersection $\mathcal F\mathcal C(F)\cap \pi^{-1}(\Sigma(M))$,
the set $\Sigma(F)$ is a closed analytic subset of $\mathcal F\mathcal C(F)$.

Denote by $\mathcal F\mathcal C^M(F)$ the flat connections whose holonomy
representations are in $\Rep^M(\Gamma;V)$.
Observe that $\mathcal F\mathcal C^M(F)$ is a closed analytic subset of 
$\mathcal F\mathcal C(F)$.

Let $g$ be a Riemannian metric on $M$ and let $\mu$ be a Hermitian fiber
metric on $F$. The Ray--Singer torsion $T_\an(\nabla,g,\mu)$ defines a 
positive real-valued function 
$$
T_\an^{g,\mu}\co\mathcal F\mathcal C(F)\to\mathbb R.
$$
The following is a straightforward
consequence of \fullref{T:ratio} and \fullref{P:hol_omega}.

\begin{corollary}\label{C:meroa}
Let $M$ be a closed connected manifold, let $F$ be a complex vector bundle over $M$,
let $g$ be a Riemannian metric on $M$ and let $\mu$ be a Hermitian fiber metric on $F$.
Then there exists a meromorphic function on $\mathcal F\mathcal C^M(F)$ 
whose zeros and poles are contained in $\Sigma(F)$ and whose restriction 
to $\mathcal F\mathcal C^M(F)\setminus \Sigma(F)$ has $T_{\an}^{g,\mu}$ as absolute value. 
\end{corollary}

\begin{proof}
Choose a base point $x_0\in M$ and an Euler structure $\e\in\Eul_{x_0}(M;\Z)$.
Choose a homology orientation $\mathfrak o$ of $M$. Let $\Gamma=\pi_1(M,x_0)$ denote
the fundamental group, and let $V=F_{x_0}$ denote the
fiber of $F$ over $x_0$. In view of 
\fullref{T:ratio}\itemref{T:ratio:i} $\mathcal T^{\e,\mathfrak o}_\comb$ is
a rational function on $\Rep^{\smash{M}}(\Gamma;V)$ with zeros and poles contained in $\Sigma(M)$. 
Since the holonomy map $\pi\co\mathcal F\mathcal C^M(F)\to\Rep^M(\Gamma;V)$ is
holomorphic, the
composition $f_1:=\mathcal T^{\e,\mathfrak o}_\comb\circ\pi$ is thus a meromorphic
function on $\mathcal F\mathcal C^M(F)$ with zeros and poles contained in $\Sigma(F)$.
Consider the co-Euler structure $\e^*:=P(\e)$ and choose $\alpha$ such that $\e^*=[g,\alpha]$.
Let $\tilde\omega\co\mathcal C(F)\to\Omega^1(M;\C)$ be an affine map as in \fullref{P:hol_omega}.
Its restriction to $\mathcal F\mathcal C(F)$ is certainly holomorphic, and thus the map
$f_2\co\mathcal F\mathcal C(F)\to\C^*$, $f_2(\nabla):=e^{S(\tilde\omega(\nabla,\mu),\alpha)}$
is holomorphic too. The product $f:=f_1\cdot f_2$ therefore is a meromorphic
function on $\mathcal F\mathcal C^M(F)$ with zeros and poles are contained $\Sigma(F)$. 
In view of \fullref{T:ratio}\itemref{T:ratio:ii} we have
$$
|f(\nabla)|=|f_1(\nabla)|e^{\Re(S(\tilde\omega(\nabla,\mu),\alpha))}
=|\mathcal T^{\e,\mathfrak o}_\comb(\rho_\nabla)|e^{S(\Re(\tilde\omega(\nabla,\mu)),\alpha)}
=|T_\an^{g,\mu}(\nabla)|
$$ 
for all $\nabla\in\mathcal F\mathcal C^M(F)\setminus\Sigma(F)$.
\end{proof}

If the dimension of $M$ is odd, the Ray--Singer torsion defines a positive 
real-valued function 
$$
T_\an\co\Rep^M(\Gamma;V)\setminus\Sigma(M)\to\mathbb R,\quad
T_\an(\rho):=T_\an(\nabla_\rho,g,\mu)
$$ 
where $g$ is any Riemannian metric on $M$, and $\mu$ is any fiber metric in $F_\rho$. 
In view of the Hermitian anomaly
\cite[Theorem~0.1]{BZ92} and because the Euler form $E(g)$ is identically
zero, this is indeed independent of $\mu$ and $g$.

\begin{corollary}\label{C:merob}
Let $M$ be an odd-dimensional closed connected manifold, let $\Gamma=\pi_1(M,x_0)$
denote its fundamental group and let $V$ be a finite-dimensional complex vector space.
Then there exists a rational function on $\Rep^M(\Gamma;V)$ whose restriction to
$\Rep^M(\Gamma;V)\setminus\Sigma(M)$ has $(T_\an)^2$ as absolute value.

If the involution $\nu$ on $\Eul(M;\Z)$ has a fixed point, or more generally, if
the canonical Euler structure $\e_\can\in\Eul(M;\R)$ is integral (\fullref{R:caneul})
then even $T_\an$ is the absolute value of a rational function on $\Rep^M(\Gamma;V)$.
\end{corollary}

\begin{proof}
Choose an Euler structure $\e\in\Eul(M;\Z)$ and a homology orientation $\mathfrak o$
of $M$. From \fullref{T:ratio}\itemref{T:ratio:i} we know that
$h_1:=(\mathcal T^{\e,\mathfrak o}_\comb)^2$ is a rational function on
$\Rep^M(\Gamma;V)$. We want to construct an additional regular function   
$h_2\co\Rep(\Gamma;V)\to\C$ so that the product $h_1\cdot h_2$ has as absolute value 
$(T_\an)^2$.

Choose a Riemannian metric $g$ on $M$. Consider the
co-Euler structure $\e^*:=P(\e)$, and choose
$\alpha\in\Omega^{n-1}(M;\mathcal O_M)$ such that $\e^*=[g,\alpha]$.
Note that $d\alpha=E(g)=0$, hence $\alpha$ defines a cohomology class 
$[\alpha]\in H^{n-1}(M;\mathcal O_M)$. We claim that the preimage $\PD^{-1}([2\alpha])\in
H_1(M;\R)$ is integral. We postpone this first step to the second half of
this proof.

Next, we choose closed loops $\gamma_i$ based at $x_0$ in $M$, $1\leq i\leq r$,
which form a set of generators of $\Gamma$. Since $\PD^{-1}([2\alpha])$ is
integral we find integers $m_i\in\Z$ such that
$\PD^{-1}([2\alpha])=\sum_{i=1}^rm_i\gamma_i$ in $H_1(M;\R)$. We then have
$$
\sum_{i=1}^rm_i\int_{\gamma_i}\omega=2\int_M\omega\wedge\alpha
$$
for every closed $1$--form $\omega\in\Omega^1(M;\R)$.
Consider the regular map
$$
h_2\co\Rep(\Gamma;V)\to\C,\quad
h_2(\rho):=\prod_{i=1}^r\bigl(\det(\rho(\gamma_i))\bigr)^{m_i}.
$$
The product $h:=h_1\cdot h_2$ is a rational function on $\Rep^M(\Gamma;V)$.
We claim that $T_\an(\rho)^2=|h(\rho)|$ for all $\rho\in\Rep^M(\Gamma;V)\setminus\Sigma(M)$.

To see this, fix $\rho\in\Rep^M(\Gamma;V)\setminus\Sigma(M)$, consider the vector
bundle $F_\rho$ over $M$ with flat connection $\nabla_\rho$ and choose a Hermitian
fiber metric $\mu$ on $F_\rho$. Recall from \fullref{R:KT} that
$$
\bigl|\det(\rho(\gamma_i))\bigr|=\exp\int_{\gamma_i}\omega(\nabla_\rho,\mu)
$$
for all $1\leq i\leq r$. We conclude
\begin{multline*}
|h_2(\rho)|=
\prod_{i=1}^r\bigl|\det(\rho(\gamma_i))\bigr|^{m_i}
=\prod_{i=1}^r\Bigl(\exp\int_{\gamma_i}\omega(\nabla_\rho,\mu)\Bigr)^{m_i}
\\=\exp\Bigl(\sum_{i=1}^rm_i\int_{\gamma_i}\omega(\nabla_\rho,\mu)\Bigr)
=\exp\Bigl(2\int_M\omega(\nabla_\rho,\mu)\wedge\alpha\Bigr)
=e^{2S(\omega(\nabla_\rho,\mu),\alpha)}.
\end{multline*}
Together with \fullref{T:ratio}\itemref{T:ratio:ii} we obtain
$$
|h(\rho)|
=|h_1(\rho)||h_2(\rho)|
=|\mathcal T^{\e,\mathfrak o}_\comb(\rho)|^2e^{2S(\omega(\nabla_\rho,\mu),\alpha)}
=T_\an(\nabla_\rho,g,\mu)^2=T_\an(\rho)^2.
$$
Thus $(T_\an)^2$ is indeed the absolute value of the rational map $h$.

Next we show that $\PD^{-1}([2\alpha])\in H_1(M;\R)$ is an integral class.
Recall from \fullref{R:caneul} that there is a canonical Euler structure
$\e_\can\in\Eul(M;\R)$ which is the unique fixed point of the involution
$\nu$ on $\Eul(M;\R)$. Note that $\e_\can$ will in general not
be an integral Euler structure. 
We identify $H_1(M;\R)$ with $\Eul(M;\R)$ by the map $a\mapsto\e_\can+a$,
and denote this identification by $\e_\can+\bullet\co H_1(M;\R)\to\Eul(M;\R)$.

Similarly we use the canonical co-Euler structure
$\e^*_\can=[g,0]\in\Eul^*(M;\R)$ of \fullref{R:cancoeul} to identify
$H^{n-1}(M;\mathcal O_M)$ with $\Eul^*(M;\R)$ by the map 
$[\beta]\mapsto\e_\can^*+[\beta]=[g,-\beta]$, and denote this identification by 
$\e^*_\can+\bullet\co H^{n-1}(M;\mathcal O_M)\to\Eul^\ast(M;\R)$.

Since from \fullref{affPD}, $P\co\Eul(M;\R)\to\Eul^*(M;\R)$ intertwines the canonical involutions
on $\Eul(M;\R)$ and $\Eul^*(M;\R)$
we must have $P(\e_\can)=\e^*_\can$.

We also use the Euler structure $\e$ to identify $H_1(M;\Z)$ with $\Eul(M;\Z)$ by the map 
$c\mapsto\e+c$, and denote this identification by $\e+\bullet\co H_1(M;\Z)\to\Eul(M;\Z)$.

Then there exists $a_0\in H_1(M;\R)$ so that the canonical map $\Eul(M;\Z)\to\Eul(M;\R)$
becomes, via the above identifications, the map $c\mapsto a_0+c$.

We summarize these observations in the following commutative diagram:
$$
\xymatrix{
\Eul(M;\Z) \ar[r] & \Eul(M;\R) \ar[r]^P & \Eul^*(M;\R)
\\
H_1(M;\Z) \ar[r]^{a_0+\bullet} \ar[u]^{\e+\bullet} 
& H_1(M;\R) \ar[r]^{\PD} \ar[u]^{\e_\can+\bullet}
& H^{n-1}(M;\mathcal O_M) \ar[u]^{\e_\can^*+\bullet}
}
$$
The commutativity of the diagram gives $\PD^{-1}([2\alpha])=-2a_0$.

To see that $2a_0$ is actually an integral homology class we note that via the
identification $\e_\can+\bullet$ the canonical involution on Euler
structures with real coefficients is given by
$\nu(a)=-a$, $a \in H_1(M;\R)$. Via the identification $\e+\bullet$  the canonical involution
on Euler structures with integer coefficients must be of the form
$\nu(c)=c_0-c$ for some $c_0\in H_1(M;\Z)$.
Since the canonical map $\Eul(M;\Z)\to\Eul(M;\R)$ intertwines these canonical involutions,
we must have $a_0+c_0=-a_0$. Hence $\PD^{-1}([2\alpha])=-2a_0=c_0$ is indeed integral.

For the last statement of the proposition, note first that any fixed point
of the involution on $\Eul(M;\Z)$ must be mapped to $\e_\can$ via
the canonical map from $\Eul(M;\Z)$ to $\Eul(M;\R)$. If the involution on $\Eul(M;\Z)$
admits a fixed point then the Euler structure $\e_\can$ must therefore be
integral. If $\e_\can$ is integral then $\PD^{-1}([\alpha])=-a_0$ must be
integral. Proceeding as above, one constructs in this case a rational
function on $\Rep^M(\Gamma;V)$ with absolute value $T_\an$.
\end{proof}

\begin{example}
Consider the circle $M=S^1$ and the one-dimensional vector space $V=\C$. 
In this case $\Gamma=\Z$, $\Rep(\Gamma;V)=\C\setminus\{0\}$,
$\Rep^M_0(\Gamma;V)=\C\setminus\{0,1\}$ and
$\Rep^M(\Gamma;V)=\C\setminus\{0\}$. From \fullref{T:marcsik}, 
or via a direct computation, we see that in this case the Ray--Singer torsion 
is given by
$$
T_\an\co\Rep^M(\Gamma;V)=\C\setminus\{0\}\to\C,\quad T_\an(z)=\frac{|z-1|}{\sqrt{|z|}}.
$$
Therefore the Ray--Singer torsion will in general not be the absolute value
of a rational function on $\Rep^M(\Gamma;V)$. The square of the Ray--Singer
torsion, however, is the
absolute value of the rational function $z\mapsto (z-1)^2/z$ which happens to
have a zero at $\Sigma(M)=\{1\}$.
\end{example}

\begin{appendix}

\section{Some homological algebra}\label{S:app}

Let $\K$ be a field of characteristics zero.
Suppose $C^*$ is a finite-dimensional $\Z$--graded complex over 
$\K$ with differential of degree one, $d^q\co C^q\to C^{q+1}$.
We will always assume $C^q=0$ for $q<0$.
Let us write $Z^q:=\ker d^q\subseteq C^q$ for the cocycles, 
$B^q:=\img d^{q-1}\subseteq C^q$ for the coboundaries, and 
$H^q:=Z^q/B^q$ for the cohomology.

Let $\cc^q$ be a base of $C^q$, let $\hh^q$ be a base of $H^q$, and let
$\bb^q$ be a base of $B^q$. Recall the short exact sequences
$$
0\to Z^q\to C^q\xrightarrow{d}B^{q+1}\to0
\quad\text{and}\quad
0\to B^q\to Z^q\to H^q\to0.
$$
Choose lifts $\tilde\hh^q$ in $Z^q$ of $\hh^q$. Then
$\bb^q\tilde\hh^q$ is a base of $Z^q$. Choose lifts $\tilde\bb^{q+1}$
in $C^q$ of $\bb^{q+1}$. Then $(\bb^q\tilde\hh^q)\tilde\bb^{q+1}$
is a base of $C^q$. For two bases $\aa_1$ and $\aa_2$ of a common
vector space we write $[\aa_1/\aa_2]$ for the determinant of the 
matrix expressing $\aa_1$ in terms of $\aa_2$. The quantity
$\prod_q[\cc^q/(\bb^q\tilde\hh^q)\tilde\bb^{q+1}]^{(-1)^q}$ does 
not depend on the choice of $\bb^q$ and does not depend on the choice
of lifts $\tilde\bb^q$ and $\tilde\hh^q$ \cite{M66}.

Recall that the determinant line of a finite-dimensional vector space 
$V$ is defined as $\det V:=\Lambda^{\dim V}V$. Moreover, we will write $L^1:=L$ 
and $L^{-1}:=\Hom(L;\K)$ for a line, \ie one-dimensional vector space, $L$.
Consider the graded determinant line of $C^*$
$$
\det C^*:=\bigotimes_q(\det C^q)^{(-1)^q}
$$
and let $\cc$ denote the base of $\det C^*$ induced from the 
bases $\cc^q$. Similarly, let $\hh$ denote the base of $\det H^*$ induced from
the bases $\hh^q$. Consider the isomorphism
\begin{equation}\label{E:rst}
\varphi_{C^*}\co \det C^*\to\det H^*,
\qquad
\cc\mapsto(-1)^N
\prod_q[\cc^q/(\bb^q\tilde\hh^q)\tilde\bb^{q+1}]^{(-1)^q}\hh
\end{equation}
where $N:=\sum_{q\geq 0}\alpha^q\beta^q$ with
$\alpha^q:=\sum_{j\geq q}\dim C^j$ and
$\beta^q:=\sum_{j\geq q}\dim H^j$ \cite{FT99,Tu01}.
The isomorphism $\varphi_{C^*}\co\det C^*\to\det H^*$ does not depend 
on the choice of bases $\cc^q$ or $\hh^q$. The sign $(-1)^N$ is referred to
as Turaev's sign refinement. Note that $(-1)^N=1$ if $C^*$ is acyclic.

Suppose we have a short exact 
sequence of finite-dimensional $\Z$--graded cochain complexes over $\K$:
\begin{equation}\label{E:ses}
0\to C_0^*\to C_1^*\to C_2^*\to0.
\end{equation} 
Then there are isomorphisms $\varphi_{C_i^*}\co\det C^*_i\to\det H^*_i$,
$i=0,1,2$. Let $\mathcal H^*$ denote the corresponding long exact sequence
\begin{equation}\label{E:lehs}
0\to H_0^0\to H_1^0\to H_2^0\xrightarrow{\delta^0}
H_0^1\to
\cdots\xrightarrow{\delta^{q-1}}
H_0^q\to H_1^q\to H_2^q\xrightarrow{\delta^q} H_0^{q+1}\to\cdots
\end{equation}
considered as (acyclic) $\Z$--graded cochain complex, the grading so that
$\mathcal H^{3q+i}=H_i^q$. As an
acyclic complex it provides an isomorphism
$\varphi_{\mathcal H^*}\co\det\mathcal H^*\to\K$ which we regard as an 
isomorphism
$\varphi_{\mathcal H^*}\co\det H_0^*\otimes\det H_2^*\to\det H^*_1$.

The short exact sequence \eqref{E:ses} provides isomorphisms
$\det C^q_0\otimes\det C^q_2=\det C^q_1$ for all $q$, and then an 
isomorphism $\psi\co\det C_0^*\otimes\det C^*_2\to\det C^*_1$.
For a proof of the following lemma consult Nicolaescu \cite[Proposition~1.18]{N03}
or look carefully at the proof of Theorem~3.2 in Milnor \cite{M66}.

\begin{lemma}\label{L:detses}
The following diagram commutes:
$$
\xymatrix{
\det C_0^*\otimes\det C_2^*
\ar[r]^-{(-1)^y\psi}\ar[d]_{\varphi_{C_0^*}\otimes\varphi_{C_2^*} }
& \det C_1^* \ar[d]^{\varphi_{C_1^*}}
\\
\det H_0^*\otimes\det H_2^* \ar[r]^-{\varphi_{\mathcal H^*}} & \det H_1^*
}
$$
Here
$
y=N(C_0^*)+N(C_1^*)+N(C_2^*)+\sum_qf_1^q\cdot b_2^q+b_2^q\cdot b_0^{q+1}+f_2^q\cdot b_0^{q+1}
$
with $b_i^q:=\dim B^q(C_i^*)$ and
$f_i^q:=\dim B^{3q+i}(\mathcal H^*)$.
\end{lemma}

Following Farber and Turaev \cite{FT99} we call the above diagram the
\emph{fusion diagram\/}, and $(-1)^y\psi$ the \emph{fusion isomorphism\/}.
Note however, that the fusion homomorphism in general depends on the
differentials of the complexes. If the short exact sequence of complexes
splits it does not \cite{FT99}.

\begin{lemma}\label{L:ee}
In the situation of \fullref{L:detses} suppose in addition that
$C_2^*$ is acyclic and two-dimensional, concentrated in degree $k$ and 
$k+1$. Then the sign $(-1)^y$ does not depend on the differentials
on $C_0^*$ and $C^*_1$.
\end{lemma}

\begin{proof}
Indeed we have 
\begin{eqnarray*}
N(C^*_0)+N(C^*_1)&=&\beta^{k+1}(C^*_0) \mod2
\\
N(C_2^*)&=&0
\\
f_1^q&=&\dim H^q_1=\dim H^q_0
\\
f_2^q&=&0
\\
b_2^q&=&0\quad\text{for $q\neq k+1$}
\\
b_2^{k+1}&=&1.
\end{eqnarray*}
Using the general fact 
$\alpha^q(C^*_0)+\beta^q(C^*_0)+b^q_0=0 \mod2$
this implies
$$
y=\beta^{k+1}(C_0^*)+\dim H^{k+1}_0+b_0^{k+2}
=\beta^{k+2}(C_0^*)+b_0^{k+2}
=\alpha^{k+2}(C_0^*)
\mod2.
$$
Hence $(-)^y$ does not depend on the differentials.
\end{proof}

Since elementary expansions give rise to short exact sequences as in 
\fullref{L:ee}, we have the following corollary.

\begin{corollary}\label{C:app}
The Milnor--Turaev torsion 
(\fullref{D:MT-tor}) does indeed not depend on the triangulation.
\end{corollary}

\end{appendix}

\bibliographystyle{gtart}
\bibliography{link}

\end{document}